\documentclass[a4paper,12pt]{article}
\usepackage[a4paper,margin=1in]{geometry}
\usepackage[utf8]{inputenc}
\usepackage[T1]{fontenc}
\usepackage[english]{babel}
\usepackage{placeins}
\usepackage{graphicx}
\usepackage{amsmath}
\usepackage{amssymb}
\usepackage{amsthm}
\usepackage[hidelinks,hypertexnames=false]{hyperref}
\usepackage{cleveref}
\usepackage{mathtools,thmtools}
\usepackage{enumerate}
\usepackage{enumitem}
\usepackage[autostyle=true,german=quotes]{csquotes}
\usepackage{color}
\usepackage{esint}
\numberwithin{equation}{section}
\crefname{subsection}{subsection}{subsections}
\crefname{subsubsection}{subsubsection}{subsubsections}

\let\originalleft\left
\let\originalright\right
\renewcommand{\left}{\mathopen{}\mathclose\bgroup\originalleft}
\renewcommand{\right}{\aftergroup\egroup\originalright}

\newcommand{\R}{{\mathbb R}}

\DeclareMathOperator{\Div}{div}

\newcommand{\opl}[1][s]{\ell({#1})}
\newcommand{\opL}[1][s]{\mathcal L({#1})}

\newcommand{\opS}[1][s]{\mathcal S({#1})}
\newcommand{\opD}[1][s]{\mathcal D({#1})}
\newcommand{\opK}[1][s]{\mathcal K({#1})}
\newcommand{\opKdual}[1][s]{\mathcal K'({#1})}
\newcommand{\opV}[1][s]{\mathcal V({#1})}
\newcommand{\opW}[1][s]{\mathcal W({#1})}
\newcommand{\CLinv}{\frac{\Re s}{\W{\wn}}}
\newcommand{\CLInv}{\frac{\W{\wn}}{\Re s}}
\renewcommand{\dim}{{n}}

\newcommand{\Interface}{\Gamma}

\newcommand{\NormHD}[4][\Gamma]%
{\left\|#4\right\|_{H^{1/2}(#1),#2,#3}}
\newcommand{\NormH}[1]%
{\left\|#1\right\|}
\newcommand{\NormHclassic}[2][\Gamma]%
{\left\|#2\right\|_{\W{H^{1/2}(#1)}}}
\newcommand{\NormHclassicDual}[2][\Gamma]%
{\left\|#2\right\|_{\W{H^{-1/2}(#1)}}}
\newcommand{\NormHalt}[1]%
{\left\lbrace\kern-2.5pt\left\vert#1\right\vert\kern-2.5pt\right\rbrace}
\newcommand{\NormHDalt}[4][\Gamma]%
{\left\lbrace\kern-2.5pt\left\vert#4\right\vert\kern-2.5pt\right\rbrace_{H^{1/2}(#1),#2,#3}}

\newcommand{\NormGD}[3][\Gamma]%
{\left\|#3\right\|_{H^{1/2}(#1),#2}}
\newcommand{\NormGN}[3][\Gamma]%
{\left\|#3\right\|_{H^{-1/2}(#1),#2}}
\newcommand{\NormHN}[4][\Gamma]%
{\left\|#4\right\|_{H^{-1/2}(#1),#2,#3}}
\newcommand{\NormHNalt}[4][\Gamma]%
{\left\lbrace\kern-2.5pt\left\vert#4\right\vert\kern-2.5pt\right\rbrace_{H^{-1/2}(#1),#2,#3}}

\newcommand{\wn}{\sigma}
\newcommand{\wnAbs}{\wn}
\newcommand{\wnlow}{\underline{\wn}}
\newcommand{\wnLow}{\overline{\wn}}
\newcommand{\gD}{g_{\mathrm{D}}}

\newcommand{\gN}{g_{\mathrm{N}}}

\renewcommand{\d}{\;\mathrm{d}}
\renewcommand{\Re}{\operatorname{Re}}
\newcommand{\conj}[1]{\overline{#1}}

\newcommand{\mean}[1]{\{\kern-4.2pt\{#1\}\kern-4.2pt\}}
\newcommand{\meanD}[2][B_R]{\mean{#2}_{\mathrm{D},#1}}
\newcommand{\meanN}[2][B_R]{\mean{#2}_{\mathrm{N},#1}}
\newcommand{\jump}[1]{\left[#1\right]}
\newcommand{\jumpD}[2][B_R]{\jump{#2}_{\mathrm{D},#1}}
\newcommand{\jumpN}[2][B_R]{\jump{#2}_{\mathrm{N},#1}}

\theoremstyle{definition}

\theoremstyle{plain}
\newtheorem{theorem}{Theorem}[section]
\newtheorem{lemma}[theorem]{Lemma}
\newtheorem{remark}[theorem]{Remark}

\newtheorem{cor}[theorem]{Corollary}

\newcommand{\Const}[1]{C_{\mathrm{#1}}}

\usepackage{totcount}
\newtotcounter{nconst}
\makeatletter
\newcount\rc@count
\let\rc@clearconstantlist\empty

\newcommand\rc@clearconstant[1]{\global\expandafter\let\csname rc@const@#1\endcsname\undefined}

\newcommand\resetconstants[1]{%
    \def\rc@constname{#1}
    \global\rc@count=1\relax 
    \bgroup 
        \let\\\rc@clearconstant 
        \rc@clearconstantlist
        \global\let\rc@clearconstantlist\empty 
    \egroup
}
\newcommand\const[1]{%
    \@ifundefined{rc@const@#1}{%
        \expandafter\xdef\csname rc@const@#1\endcsname{%
           \noexpand\rc@useconst{\rc@constname}{\the\rc@count}%
        }%
        \g@addto@macro\rc@clearconstantlist{\\{#1}}%
		\global\setcounter{nconst}{\the\rc@count}\relax
        \global\advance\rc@count1\relax
    }{}%
    \csname rc@const@#1\endcsname
}

\newcommand\rc@useconst[2]{\ensuremath{#1_{#2}}}
\makeatother

	\resetconstants{C}

\newcommand{\W}[1]{{#1}}

\title{Optimal trace norms for Helmholtz problems}
\author{Benedikt Gr\"a{\ss}le\footnote{Institut für Mathematik, Universität Zürich, Winterthurerstr.~190, CH-8057 Zürich, Switzerland \\
  (benedikt.graessle@math.uzh.ch)}
}
\date{}%
\begin{document}

\maketitle

\begin{abstract}
	The natural $H^1(\Omega)$ energy norm for Helmholtz problems is weighted 
	with the wavenumber modulus $\wn$
	and induces \W{weighted norms} on the trace spaces $H^{\pm1/2}(\Gamma)$ 
	by minimal extension to $\Omega\subset\R^n$.
	This paper \W{provides} an explicit characterisation \W{through}
	weighted Sobolev-Slobodeckij norms and scaling estimates,
	highlighting \W{the} dependence on the geometry of the extension set 
	$\Omega\subset\R^n$ and the weight $\wn$.
	The analysis identifies conditions under which these trace norms are 
	intrinsic to the isolated boundary component $\Gamma\subset\partial\Omega$ 
	and \W{establishes $\wn$-explicit trace estimates} in weighted spaces.
	\W{In these norms, the Helmholtz potential and boundary integral operators satisfy
	improved coercivity and continuity estimates without additional low-frequency
	factors that deteriorate as $\wn\to 0$; for $n\geq 3$, the same analysis also improves
	the corresponding bounds in the classical unweighted trace norms.
}
\end{abstract}



\section{Introduction}\label{sec:Introduction}
Wave propagation phenomena
appear across various scientific and engineering disciplines, 
from acoustics and electromagnetism to elasticity and quantum mechanics.
A fundamental model for time-harmonic waves is the Helmholtz equation 
for a complex wavenumber $s$,
which encodes both the frequency of oscillations 
as well as damping or absorption effects from the ambient medium.
The mathematical analysis and numerical theory of Helmholtz-type problems 
rely on variational formulations in appropriate Sobolev spaces equipped with
weighted norms that reflect the oscillatory nature of waves.

Time-harmonic waves in an open set $\Omega\subset\R^n$ are described by the Helmholtz operator
\begin{align}\label{eqn:HH_op}
	v\mapsto \opL v\coloneqq -\Delta v + s^2 v,
	\qquad s\in \mathbb{C}\setminus\{0\}
\end{align}
and possibly boundary conditions (e.g., of Dirichlet, Neumann, or impedance type).
The \emph{natural energy norms} for $H^1(\Omega)$ and $H(\Omega,\Div)$ associated to the Helmholtz 
operator~\eqref{eqn:HH_op} are weighted with the wavenumber modulus $\wn=|s|>0$~as
\begin{align}\label{eqn:H1_norm}
	\|v\|_{H^1(\Omega),\wn}^2 
	&= \|\nabla v\|_{L^2(\Omega)}^2 + 
	\wn^2\|v\|_{L^2(\Omega)}^2
	&&\text{for all }v\in H^1(\Omega),\\
	\|q\|_{H(\Omega,\Div),\wn}^2 
	&= \wn^{-2}\|\Div q\|_{L^2(\Omega)}^2 + 
	\|q\|_{L^2(\Omega)}^2
	&&\text{for all }q\in H(\Omega,\Div).\label{eqn:Hdiv_norm}
\end{align}
\W{In fact,~\eqref{eqn:H1_norm} is the coercive energy norm associated with the sesquilinear form of~\eqref{eqn:HH_op},
and~\eqref{eqn:Hdiv_norm} the corresponding flux energy: If $\opL v=0$, then $\|\nabla
v\|_{H(\Omega,\Div),\wn}=\|v\|_{H^1(\Omega),\wn}$.
The weight in \eqref{eqn:H1_norm} balances the $L^2$-term with the gradient
term at wavenumber~$s$.

The analysis of boundary conditions, 
or interface conditions} and the coupling with exterior problems,
requires to understand the Dirichlet trace space $H^{1/2}(\Gamma)$ and the Neumann trace space 
$H^{-1/2}(\Gamma)$ on boundary parts $\Gamma\subset\partial\Omega$.
These spaces describe the admissible boundary data for solutions (and their fluxes) and are central 
in variational formulations%
~\cite{Gri:EllipticProblemsNonsmooth1985,BH:FormulationVariationnelleEspacetemps1986,Ned:AcousticElectromagneticEquations2001}.
They also arise naturally in the study of boundary integral equations%
~\cite{Cos:BoundaryIntegralOperators1988,Von:BoundaryIntegralEquations1989,McL:StronglyEllipticSystems2000,CM:WavenumberexplicitBoundsTimeharmonic2008,CK:InverseAcousticElectromagnetic2019}, 
which represent solutions by its Cauchy traces on $\Gamma$ and
enable fast numerical solvers based on a discretisation of the compact lower-dimensional boundary%
~\cite{SS:BoundaryElementMethods2011,LM:WavenumberexplicitHpBEMHigh2011,CH:IntegralEquationsMultiScreens2013}.
\W{A central focus in the standard analysis}
(see~\cite{CM:WavenumberexplicitBoundsTimeharmonic2008,LM:WavenumberexplicitHpBEMHigh2011,MS:WavenumberExplicitConvergence2011,SW:WavenumberexplicitParametricHolomorphy2023})
is on wavenumber-explicit estimates to characterise the dependence of important operators on the wavenumber.
However, the estimates involving traces often depend explicitly on a lower bound of $\sigma=|s|$, namely
\begin{align}\label{eqn:wnlow_def}
	\wnlow\coloneqq\min\{1,\wn\}
	\qquad\text{or its inverse}\quad
	\wnLow\coloneqq\wnlow^{-1}=\max\{1,\wn^{-1}\},
\end{align} 
and deteriorate for $\sigma\to0$
(or are only stated for sufficiently large wavenumbers $\sigma\geq\sigma_0$), cf.%
~\cite[sec.~4.4]{BS:IntegralEquationMethods2022} for references in the context of boundary integral
operators.
Such estimates can be sharpened and the dependence on~\eqref{eqn:wnlow_def} removed by
endowing $H^{\pm1/2}(\Gamma)$ with natural weighted trace norms (corresponding to~\eqref{eqn:H1_norm}--\eqref{eqn:Hdiv_norm})
introduced and analysed in this paper.
Weighted trace norms have already been employed for a sharp wavenumber-explicit analysis, e.g.,
in~\cite{MS:WavenumberExplicitConvergence2011} (and for the related Maxwell problem
in~\cite{MS:WavenumberExplicitHpFEMAnalysis2024}),
but their weight is \emph{not} optimal for all $\sigma>0$ shown by the assumption $\sigma\geq\sigma_0>0$ therein.

This paper contributes to the mathematical analysis of Helmholtz PDE by providing
$\wn$-explicit characterisations and estimates for the natural trace norms induced by the Helmholtz operator.
The results enable sharper wavenumber-explicit estimates 
in optimally weighted Sobolev norms~\cite{GHS:StableSkeletonIntegral2025}
that may be used to improve existing results in the literature.
As an application, the paper presents novel 
operator bounds for boundary integral operators 
that are independent of $\wnlow,\wnLow$ in~\eqref{eqn:wnlow_def} 
and valid for all $\sigma>0$.
While the focus is on the scalar Helmholtz problem, the techniques and results extend to vector-valued
function spaces and the related Maxwell setting.

\medskip
\noindent\textbf{Natural trace norms.}
For smooth (Lipschitz) boundaries $\Gamma\subset\partial\Omega$, $H^{1/2}(\Gamma)$ 
consists of all traces of $H^{1}(\Omega)$ functions on $\Gamma$.
This identification implies the natural isomorphy
\begin{align}\label{eqn:H12_quotient space}
	H^{1/2}(\Gamma)\cong H^1(\Omega)/H^1_{0,\Gamma}(\Omega)
\end{align}
with the quotient space over the kernel $H^1_{0,\Gamma}(\Omega)$ of the Dirichlet trace map.
The quotient space inherits its norm from~\eqref{eqn:H1_norm} and the identification%
~\eqref{eqn:H12_quotient space} induces a natural trace norm in 
$H^{1/2}(\Gamma)$ by minimal extension
\begin{align}\label{eqn:H12_norm_def}
	\NormHD[\Gamma]{\Omega}{\wn}\gD
	\coloneqq\inf_{\substack{v\in H^1(\Omega)\\ v|_{\Gamma} = \gD}}
		\|v\|_{H^1(\Omega),\wn}
	\qquad\text{for all }\gD\in H^{1/2}(\Gamma).
\end{align}
The relation~\eqref{eqn:H12_quotient space} may also be employed to define
$H^{1/2}(\Gamma)$ with a weighted norm~\eqref{eqn:H12_norm_def} for general (rough) open 
$\Omega\subset\R^n$ that do
not admit classical traces%
~\cite{CH:IntegralEquationsMultiScreens2013,CH:IntegralEquationsElectromagnetic2016,HPS:TracesHilbertComplexes2023}.

To simplify the discussion, it is assumed throughout this paper that $\Gamma\subset\partial\Omega$ is
an \emph{isolated boundary component} of some Lipschitz set $\Omega\subset\R^n$ (called extension set)
in the sense that it is nonempty and \W{the distance
$\mathrm{dist}(\Gamma,\partial\Omega\setminus\Gamma)>0$ is positive whenever $\partial\Omega\setminus\Gamma\ne\emptyset$.
(For a simpler first reading, 
	the reader may keep in mind the prototypical interface $\Gamma=\partial\Omega$
with $\Omega$ bounded Lipschitz and connected exterior $\R^n\setminus\Omega$.)
In this setting,} $H^{-1/2}(\Gamma)=(H^{1/2}(\Gamma))'$ is the dual space of $H^{1/2}(\Gamma)$,
associated to the normal traces of $H(\Omega,\Div)$ vector fields
and naturally equipped with the operator norm
\begin{align}\label{eqn:H12_dual_norm_def}
	\NormHN[\Gamma]{\Omega}{\wn}{\gN}
	&\coloneqq
	\sup_{\substack{0\ne \gD\in H^{1/2}(\Gamma)}}
	\frac{\left|\left\langle \gN,\gD\right\rangle_{\Gamma}\right|}
	{\NormHD[\Gamma]{\Omega}{\wn}{\gD}}
	\qquad\text{for all }\gN\in H^{-1/2}(\Gamma),
\end{align}
where $\left\langle \bullet,\bullet\right\rangle_{\Gamma} $ denotes the dual pairing in
$H^{-1/2}(\Gamma)\times H^{1/2}(\Gamma)$.
Standard variational principles identify~\eqref{eqn:H12_dual_norm_def}
as a minimal extension norm for the extension by $H(\Omega,\Div)$ vector-fields
equipped with the weighted norm~\eqref{eqn:Hdiv_norm}, subject to boundary conditions  on
$\partial\Omega\setminus\Gamma$.

The weighted trace norms~\eqref{eqn:H12_norm_def}--\eqref{eqn:H12_dual_norm_def} are
naturally induced by the energy norms~\eqref{eqn:H1_norm}--\eqref{eqn:Hdiv_norm} and
encode the correct scaling of the solutions to the Helmholtz equation in $\wn$.
For $\wn=1$, these are standard norms for $H^{\pm1/2}(\Gamma)$ 
on Lipschitz boundaries $\Gamma\subset\partial\Omega$ and
equivalent, e.g., to classical 
interpolation or intrinsic Sobolev-Slobodeckij norms%
~\cite{LM:NonhomogeneousBoundaryValue1972,McL:StronglyEllipticSystems2000,BF:MathematicalToolsStudy2013}.
However,
special care is required in the mathematical analysis of these norms for general $\wn>0$
as their definitions~\eqref{eqn:H12_norm_def}--\eqref{eqn:H12_dual_norm_def} are \emph{extrinsic}:
The trace norms depend on the particular extension set $\Omega$.
Different extension sets lead to norms that are in general \emph{not} $\wn$-equivalent 
in the sense that the equivalence constants may depend on $\wn$ and are not uniformly bounded for all $\wn>0$.

\medskip\noindent
\textbf{Main results.}
This paper presents a detailed $\wn$-explicit analysis of the minimal extension
norms~\eqref{eqn:H12_norm_def}--\eqref{eqn:H12_dual_norm_def} to understand the dependence on
the extension sets and the correct scaling in the weight $\wn$ induced in the trace spaces $H^{\pm1/2}(\Gamma)$.

\begin{enumerate}[label=(\alph*),align=left,leftmargin=*]
	\item 
The main result is a $\wn$-explicit comparison of~\eqref{eqn:H12_norm_def} to
a weighted intrinsic Sobolev-Slobodeckij-norm for $H^{1/2}(\Gamma)$ given by
\begin{align}\label{eqn:Sob_norm_def}
	\NormGD[\Gamma]{\wn}\gD\coloneqq\sqrt{|\gD|_{H^{1/2}(\Gamma)}^2 + \wn\min\{1,\wn\}\|\gD\|_{L^2(\Gamma)}^2}
	\qquad\text{for all }\gD\in H^{1/2}(\Gamma)
\end{align}
with the classical Gagliardo seminorm $|\bullet|_{H^{1/2}(\Gamma)}$.
\W{The scaling of~\eqref{eqn:H12_norm_def} in the weight $\wn$ depends on the particular extension 
set $\Omega$ of $\Gamma$. This}
is highlighted by the $\wn$-explicit characterisation of~\eqref{eqn:H12_norm_def}:
For bounded Lipschitz domains $\Omega\subset\R^n$ and $n\geq2$, it~holds
\begin{align*}%
		\NormHD[\Gamma]{\Omega}{\wn}\gD
		\approx
		\NormGD[\Gamma]{\wn}\gD
	\end{align*}
and, for exterior Lipschitz domains $\Omega\subset\R^n$ (with $\R^n\setminus\Omega$ bounded) and $n\geq3$,
\begin{align*}%
	\NormHD[\Gamma]{\Omega}{\wn}\gD
	\approx
	\NormGD[\Gamma]{\max\{1,\wn\}}\gD.
\end{align*}
The notation \enquote{$\approx$} entails \W{equivalence constants independent of $\sigma>0$, but possibly depending}
on $\Gamma$ and $\Omega$.
For general (Lipschitz) extension sets $\Omega$ of $\Gamma$,
the corresponding trace norms are always $\wn$-equivalent up to a multiplicative factor of%
~\eqref{eqn:wnlow_def}.
In other words, 
the weight induced by the norms~\eqref{eqn:H12_norm_def}--\eqref{eqn:H12_dual_norm_def}
in $H^{\pm1/2}(\Gamma)$ is intrinsic to $\Gamma$ 
up to~\eqref{eqn:wnlow_def}.
\W{The $\wn$-equivalence holds even in a strict sense 
if one restricts to traces from either \emph{i)}
only bounded domains, or \emph{ii)} only exterior domains (and
$n\geq3$).}

	\item The characterisation in (a) provides further insights into the weighted trace norms.
		To relate these norms for general extension sets $\Omega$ with 
		the \W{(unweighted) classical norms 
		$\NormH{\bullet}_{H^{\pm1/2}(\Gamma)}$,}
		it is essential to understand their scaling with respect to
		$\wn$ given~by
		\W{\begin{equation}
				\label{eqn:scaling_HN_bounded}
			\begin{aligned}
				\wnlow\NormH{\bullet}_{H^{1/2}(\Gamma)}
				&\lesssim
				\NormHD[\Gamma]{\Omega}{\wn}\bullet
				\lesssim \max\{1,\wn^{1/2}\}\NormH{\bullet}_{H^{1/2}(\Gamma)},\\
				\wnlow\NormHN[\Gamma]{\Omega}{\wn}\bullet
				&\lesssim
				\NormH{\bullet}_{H^{-1/2}(\Gamma)}
				\lesssim\max\{1,\wn^{1/2}\}\NormHN[\Gamma]{\Omega}{\wn}\bullet.
			\end{aligned}
		\end{equation}
		For exterior domains $\Omega$ and $n\geq3$, one obtains an improved scaling 
		without $\wnlow$:
		\begin{equation} \label{eqn:scaling_HN_bounded_improved}
			\begin{aligned}
				\NormH{\bullet}_{H^{1/2}(\Gamma)}
				&\lesssim
				\NormHD[\Gamma]{\Omega}{\wn}\bullet
				\lesssim \max\{1,\wn^{1/2}\}\NormH{\bullet}_{H^{1/2}(\Gamma)},\\
				\NormHN[\Gamma]{\Omega}{\wn}\bullet
				&\lesssim
				\NormH{\bullet}_{H^{-1/2}(\Gamma)}
				\lesssim\max\{1,\wn^{1/2}\}\NormHN[\Gamma]{\Omega}{\wn}\bullet.
			\end{aligned}
		\end{equation}%
		}%
		Trace estimates control the 
		Dirichlet traces of $H^1$ functions and the normal traces of $H(\Div)$ vector-fields
		defined on open sets $\widetilde{\Omega}\subset\R^n$ (possibly) different from the extension set
		in the natural energy norms~\eqref{eqn:H1_norm}--\eqref{eqn:Hdiv_norm}.
		Particular care is taken to identify conditions on $\widetilde{\Omega}$ and the extension set $\Omega$
		for $\wn$-independent trace constants.
		In general, the optimal constants in the trace estimate depend on~\eqref{eqn:wnlow_def}\W{,
	but there are important configurations that allow $\sigma$-independent trace estimates.}

	\item 
		The well-known single and double layer potentials associated to~\eqref{eqn:HH_op} for 
		$\Re s>0$,
		\begin{align}\label{eqn:opSD}
			\opS:H^{-1/2}(\Gamma)\to H^1(\R^n\setminus\Gamma)
			\qquad\text{and}\qquad
			\opD:H^{1/2}(\Gamma)\to H^1(\R^n\setminus\Gamma),
		\end{align}
		are the fundamental solution operators for transmission problems in $\R^n\setminus\Gamma$ across some 
		Lipschitz interface $\Gamma=\partial\Omega$%
		~\cite{McL:StronglyEllipticSystems2000,SS:BoundaryElementMethods2011,BS:IntegralEquationMethods2022}.
		Their (averaged) traces define the boundary layer operators~\cite[eqn.~(3.6)]{SS:BoundaryElementMethods2011}
		on $H^{\pm1/2}(\Gamma)$,
		often denoted by $\opV$, $\opK$, $\opKdual$, and $\opW$.
		With respect to the optimally weighted norms discussed in this paper, 
		these operators are bounded exclusively in terms of the coercivity constant $\Re s/\W{\wn}$
		of the Helmholtz operator and otherwise independent of $\sigma>0$.
		\W{As a result, these operators remain uniformly bounded as $\wn=|s|\to 0$
		along sequences with uniformly positive $\Re s/\W{\wn}$,
		and this even holds with respect to the (unweighted) classical norms if 
		$\Gamma$ is the boundary of an exterior domain and $n\geq3$.}

		These bounds are extended in~\cite[sec.~4]{GHS:StableSkeletonIntegral2025} (based on the results herein)
		to the indefinite case $\Re s=0$ and (possibly) heterogeneous coefficients.
\end{enumerate}

\medskip
\noindent\textbf{Outline.}
\Cref{sub:General notation} introduces the geometric setting and classical results on
Sobolev spaces and their traces.

\W{
\Cref{sub:Charaterisation of minimal extension norms} characterises the
weighted trace and operator norms in \cref{lem:minimal extension norm}. 
In the case $\Gamma\neq \partial\Omega$, it also introduces
an alternative trace norm 
in which the minimal extension vanishes on $\Gamma^c=\partial\Omega\setminus\Gamma$. 
The associated operator norm then relates to a minimal extension over all 
$H(\Omega,\Div)$ vector fields,
without imposing boundary conditions on $\Gamma^c$. This complements
\eqref{eqn:H12_norm_def}, where the trace on $\Gamma^c$ is left free and the
minimal flux extension associated with the operator norm
\eqref{eqn:H12_dual_norm_def} satisfies homogeneous normal conditions on $\Gamma^c$.}
\Cref{lem:Fourier_halfspace} states an explicit representation of~\eqref{eqn:H12_norm_def} on the half-space
$\R^{n}_+$ in terms of the Fourier transform\W{; see}, e.g., 
from~\cite{BH:FormulationVariationnelleEspacetemps1986,Say:RetardedPotentialsTime2016}.
The main result, \cref{thm:characterisation_intrinsic}, compares~\eqref{eqn:H12_norm_def} with 
the weighted Sobolev-Slobodeckij norm~\eqref{eqn:Sob_norm_def} 
for connected extension sets $\Omega$ 
and \W{prepares the $\wn$-explicit characterisation for extension domains
in \cref{thm:characterisation}; see part~(a).}
\Cref{rem:multiply_connected} discusses an immediate generalisation to \W{multiply connected}
extension sets.

\W{\Cref{sec:Properties of minimal extension norms} 
	establishes properties of the minimal extension norms
	required in the mathematical analysis of Helmholtz problems and departs 
	in \cref{lem:scaling in the weight,cor:scaling in the weight}
	with their scaling in the weight.
The explicit dependence of~\eqref{eqn:H12_norm_def} on the extension set $\Omega$ 
is analysed in \cref{thm:extension_sets} and 
enables the $\wn$-explicit trace control in \cref{lem:trace_inequality}; see part~(b).
Particular attention is paid to identify geometric conditions 
on the extension sets $\Omega$ that ensure $\wn$-independent trace constants.}

The paper concludes in \cref{sub:Operator bounds for boundary integral operators} 
with $\wn$-explicit continuity estimates 
for boundary integral operators in \W{weighted trace norms
and (for $n\geq3$) with improved bounds in the classical unweighted trace norms; see part~(c).
These estimate refine the known bounds from%
~\cite{BH:FormulationVariationnelleEspacetemps1986,BS:IntegralEquationMethods2022,FHS:SkeletonIntegralEquations2024},
which are recovered from part~(b) through the $\wn$-explicit scaling of the weighted trace norms.}

\section{Preliminaries}%
\label{sub:General notation}
Throughout this paper, $\wn>0$ denotes a positive weight representing the wavenumber modulus.
Standard notation for (complex-valued) Lebesgue and Sobolev spaces and 
their norms applies in $n\geq2$ dimensions. 
A domain is an open and connected set $\omega\subset\R^n$. 
It is called an exterior domain if additionally $\R^n\setminus\omega$ is bounded.
A Lipschitz set is an open set $\Omega\subset\R^n$
with compact boundary $\partial\Omega$ 
that is locally the graph of a Lipschitz function~\cite[def.~3.28]{McL:StronglyEllipticSystems2000}.
The complex-conjugation of $z\in\mathbb{C}$ is denoted by $\conj z$.

\W{Let $\emptyset\ne\Gamma\subset\partial\Omega$ be an isolated boundary component 
of $\partial\Omega$
in the aforementioned sense,
i.e., $\operatorname{dist}(\Gamma,\partial\Omega\setminus\Gamma)>0$
if $\Gamma\ne\partial\Omega$.}
In particular, both $\Gamma$ and $\partial\Omega\setminus\Gamma$ are (relatively) open and closed.
The surjective (Dirichlet) trace operator 
$(\bullet)|_{\Gamma}: H^1(\Omega)\to H^{1/2}(\Gamma)$
is the unique continuous linear extension of the pointwise
trace of smooth functions. Hence 
\begin{align}\label{eqn:C_tr_1}
	\NormHclassic[\Gamma]{v|_{\Gamma}}
	\leq \Const{tr,1}\|v\|_{H^1(\Omega),1}
	\qquad\text{for all }v\in H^1(\Omega).
\end{align}
Here $\NormHclassic[\Gamma]\bullet^2\coloneqq|\bullet|_{H^{1/2}(\Gamma)}^2 +
\|\bullet\|_{L^2(\Gamma)}^2\W{\equiv\NormGD[\Gamma]{1}\bullet^2}$
denotes the classical (unweighted) Sobolev-Slobodeckij norm for $H^{1/2}(\Gamma)$ 
\W{that is} defined 
in terms of the Gagliardo seminorm 
\begin{align}\label{eqn:Gagliardo_def}
	|\gD|_{H^{1/2}(\Gamma)}^2
	\coloneqq
	\int_{\Gamma}\int_{\Gamma}\frac{|\gD(x)-\gD(y)|^2}{|x-y|^n}\d y\d x
	\qquad\text{for all }\gD\in H^{1/2}(\Gamma).
\end{align}
\W{Its associated operator norm is denoted by $\NormH{\bullet}_{H^{-1/2}(\Gamma)}$.}
The kernel of $|\bullet|_{H^{1/2}(\Gamma)}$ coincides with the \W{constant functions over $\Gamma$
(also in the case of multiply connected $\Gamma$)}.
The compact embedding
$H^{1/2}(\Gamma)\hookrightarrow L^2(\Gamma)$
leads to a Poincar\'e inequality on $\Gamma$, namely
\begin{align}\label{eqn:Poincare_Gamma_def}
	\left\|{g}-\frac{1}{|\Gamma|}\int_\Gamma g\d s\right\|_{L^2(\Gamma)}
	\leq\Const{P,\Gamma}|{g}|_{H^{1/2}(\Gamma)}
	\qquad\text{for all }g\in H^{1/2}(\Gamma).
\end{align}
For a partition of $\Gamma=\cup_{j=1}^J\Gamma_j$ into $J\in\mathbb{N}$ disjoint open parts $\Gamma_j\subset\Gamma$ 
with positive $\operatorname{dist}(\Gamma_j,\Gamma\setminus\Gamma_j)>0$ for $j=1,\dots,J$,
there exists
a constant $\Const{spl}>0$ with%
\footnote{A proof can be found, e.g., in%
~\cite[pp.~525-526]{SS:BoundaryElementMethods2011}.}
\begin{align}\label{eqn:Gagliardo_split}
	\sum^{J}_{j=1} |\gD|_{H^{1/2}(\Gamma_j)}^2 
	\leq
	|\gD|_{H^{1/2}(\Gamma)}^2
	\leq
	\sum^{J}_{j=1} |\gD|_{H^{1/2}(\Gamma_j)}^2 
	+ \Const{spl}\|\gD\|_{L^2(\Gamma)}^2.
\end{align}
The trace operator admits a continuous right-inverse $R_\Gamma:H^{1/2}(\Gamma)\to H^1(\Omega)$ with
\begin{align}\label{eqn:R_Gamma_def}
	(R_{\Gamma}g)|_{\Gamma} = g
	\quad\text{and}\quad
	\|R_\Gamma g\|_{H^1(\Omega),1}
	\leq \Const{R,1}\NormHclassic[\Gamma]{g}
	\qquad\text{for all }g\in H^{1/2}(\Gamma).
\end{align}

Given an oriented unit normal $\nu$ on $\Gamma$, 
the normal (Neumann) trace $(\bullet\cdot\nu)|_{\Gamma}:H(\Omega,\Div)\to H^{-1/2}(\Gamma)$
maps a vector field 
$q\in H(\Omega,\Div)
\coloneqq\{q\in L^2(\Omega;\R^\dim)\ :\ \Div q\in L^2(\Omega)\}$
with square-integrable divergence 
to its normal trace $(q\cdot\nu)|_{\Gamma}$ 
and is surjective onto
the topological dual space $H^{-1/2}(\Gamma)=(H^{1/2}(\Gamma))'$.
The bilinear version of the $L^2(\Gamma)$-scalar product\W{,} given by
\begin{align*}
	\left\langle v,w\right\rangle_{\Gamma}
	\coloneqq \int_\Gamma v\, w\d s
	\qquad\text{for all }v,w\in L^2(\Gamma),
\end{align*}
extends to the bilinear dual pairing $H^{-1/2}(\Gamma)\times H^{1/2}(\Gamma)$ with the same notation.
The subspaces with boundary conditions on \W{the} isolated boundary component 
$\gamma\subset\partial\Omega$ \W{read}
\begin{align*}
	H^1_{0,\gamma}(\Omega)&\coloneqq\{v\in H^1(\Omega)\ :\ v|_{\gamma}\equiv 0 \},\\
	H_{0,\gamma}(\Omega,\Div)&\coloneqq\{q\in H(\Omega,\Div)\ :\ (q\cdot\nu)|_{\gamma}\equiv 0
	\}.
\end{align*}
Gauss' theorem reveals a duality of the boundary conditions in the sense that
\begin{align}\label{eqn:Gauss_dualtiy}
	\left\langle (q\cdot\nu)|_{\Gamma}, v|_{\Gamma}\right\rangle 
	= \int_{\Omega} v\Div q + q\cdot\nabla v\d x
	\qquad\text{for all }v\in H^1_{0,\gamma}(\Omega),q\in H_{0,\partial\Omega\setminus\gamma}(\Omega,\Div).
\end{align}
\W{An isolated boundary component $\gamma\subset\partial\Omega$ of the Lipschitz set $\Omega$
is nonempty if and only if $|\gamma|>0$.}
The Poincar\'e and, for $\gamma\ne\emptyset$, the Friedrichs inequality for bounded $\Omega$ read
\begin{align}\label{eqn:C_P}
	\left\|v-\frac{1}{|\Omega|}\int_{\Omega}v\d x\right\|_{L^2(\Omega)}
	&\leq\Const{P,\Omega}\|\nabla v\|_{L^2(\Omega)}
	&&\text{for all }v\in H^1(\Omega),\\\label{eqn:C_F}
	\|v\|_{L^2(\Omega)}
	&\leq\Const{F}\|\nabla v\|_{L^2(\Omega)}
	&&\text{for all }v\in H^1_{0,\gamma}(\Omega)
\end{align}

The Dirichlet jumps and averages of $v\in H^1(\R^n\setminus\Gamma)$ across $\Gamma\subset\partial\Omega$
are defined by
\begin{equation}\label{eqn:jump_mean_def}
	\begin{aligned}
		\jumpD[\Gamma]{v}
		&\coloneqq v^{0}|_{\Gamma} - v^{\mathrm{ext}}|_{\Gamma},
		\quad
		&\meanD[\Gamma]{v}
		&\coloneqq \tfrac{1}{2}(v^{0}|_{\Gamma} + v^{\mathrm{ext}}|_{\Gamma})
	\end{aligned}
\end{equation}
for $v^{0}\coloneqq v|_{\Omega}$ and $v^{\mathrm{ext}}\coloneqq v|_{\R^n\setminus\overline{\Omega}}$.
If additionally $\Delta v\in L^2(\R^n\setminus\Gamma)$,
the Neumann jumps and averages are given with the normal derivative $\partial_{\nu}\bullet\coloneqq\nabla
\bullet\cdot\nu$ by
\begin{align}\label{eqn:N_jump_average}
		\jumpN[\Gamma]{v}
		&\coloneqq\partial_{\nu} v^0|_{\Gamma} - \partial_{\nu}v^{\mathrm{ext}}|_{\Gamma},
		\quad
		&\meanN[\Gamma]{v}
		&\coloneqq \tfrac{1}{2}(\partial_{\nu} v^0|_{\Gamma} +
		\partial_{\nu}v^{\mathrm{ext}}|_{\Gamma}).
\end{align}

All constants $\Const{}$ are independent of $\wn$ and the notation $A\lesssim B$ 
hides a generic constant $C$ in the estimate $A\leq C B$, whereas $A\approx B$ abbreviates $A\lesssim B\lesssim A$.
The context-sensitive notation \W{$|\bullet|$ may refer to 
the $m$-dimensional (Lebesgue or surface) measure $|\omega|$ of 
some $m$-dimensional manifold $\omega\subset\R^n$ for $m\in\mathbb{N}_0$
or the absolute value $|z|$ of a complex number $z\in\mathbb{C}$}.
The $m$-dimension ball with radius $R>0$ is denoted by $B_R^m=\{x\in \R^m\ :\ \|x\|\leq R\}$ and in the case $m=n$
occasionally written $B_R=B_R^n$.

\begin{remark}[operator norm estimates]\label{rem:operator_norm_est}
	Let $\|\bullet\|_A$ and $\|\bullet\|_B$ denote two norms for $H^{1/2}(\Gamma)$ and
	$\|\bullet\|_{A,*}$ and $\|\bullet\|_{B,*}$ the associated (dual) operator norms for $H^{-1/2}(\Gamma)$.
	Any estimate of the form
	\begin{align*}
		\|\gD\|_{A}
		\leq \Const{}(A,B)\|\gD\|_{B}
		\qquad\text{for all }\gD\in H^{1/2}(\Gamma)
	\end{align*}
	implies, by definition of the operator norm, the reverse estimate 
	\begin{align*}
		\|\gN\|_{B,*}
		\leq \Const{}(A,B)\|\gN\|_{A,*}
		\qquad\text{for all }\gN\in H^{-1/2}(\Gamma)
	\end{align*}
	for the operator norm and vice versa.
\end{remark}

\section{Representation and intrinsic characterisation}%
\label{sub:Charaterisation of minimal extension norms}
This section characterises and represents the weighted trace norms%
~\eqref{eqn:H12_norm_def}--\eqref{eqn:H12_dual_norm_def}
on the isolated boundary component $\Gamma\subset\partial\Omega$ of some
Lipschitz set $\Omega\subset\R^n$ for $n\geq2$ with oriented normal $\nu$
in terms of the weighted Sobolev-Slobodeckij norm~\eqref{eqn:Sob_norm_def}.
\W{A (possibly)} nonempty boundary complement 
$\Gamma^c\coloneqq\partial\Omega\setminus\Gamma$ and~\eqref{eqn:Gauss_dualtiy}
motivate an alternative minimal extension norm for $H^{1/2}(\Gamma)$ and corresponding operator norm
for $H^{-1/2}(\Gamma)$ by
\begin{align}\label{eqn:H12_alt_def}
	\NormHDalt[\Gamma]{\Omega}{\wn}{\gD}
	&\coloneqq
	\inf_{\substack{v\in H^1(\Omega)\\ 
			v|_{\Gamma}=\gD\text{ on }\Gamma\\
	v|_{\Gamma^c}=0\text{ on }\Gamma^c}}
	\|v\|_{H^1(\Omega),\wn}
	\qquad\text{for all }\gD\in H^{1/2}(\Gamma),\\
	\NormHNalt[\Gamma]{\Omega}{\wn}{\gN}
	&\coloneqq
	\sup_{\substack{0\ne \gD\in H^{1/2}(\Gamma)}}
	\frac{\left|\left\langle \gN,\gD\right\rangle_{\Gamma}\right|}
	{\NormHDalt[\Gamma]{\Omega}{\wn}{\gD}}\label{eqn:H12_dual_alt_def}
	\qquad\text{for all }\gN\in H^{-1/2}(\Gamma).
\end{align}
\W{The two norms~\eqref{eqn:H12_norm_def} and~\eqref{eqn:H12_alt_def}
	encode two complementary boundary conditions on $\Gamma^c$:
	in~\eqref{eqn:H12_norm_def}, the trace on $\Gamma^c$ is free and
	the minimal extension satisfies homogeneous Neumann conditions on $\Gamma^c$,
	while~\eqref{eqn:H12_alt_def} enforces homogeneous Dirichlet conditions.
	The corresponding operator norms reverse this behaviour:
	the first imposes a homogeneous normal conditions on $\Gamma^c$, while the second leaves the flux unrestricted on
	$\Gamma^c$.
}

\begin{theorem}[minimal extension norm]\label{lem:minimal extension norm}
	For any $\gD\in H^{1/2}(\Gamma)$ and $\sigma>0$, the infima in~\eqref{eqn:H12_norm_def}
	and~\eqref{eqn:H12_alt_def} are attained
	for the unique solution $u\in H^1(\Omega)$ to
	\begin{align}\label{eqn:min_PDE}
		-\Delta u + \wn^2 u&=0\quad\text{in }\Omega,
						   &
		u|_{\Gamma}&=\gD\quad\text{on }\Gamma,
				   &
		(\nabla u\cdot\nu)|_{\Gamma^c} &= 0\quad\text{on }\Gamma^c;\\
		-\Delta u + \wn^2 u&=0\quad\text{in }\Omega,
						   &
		u|_{\Gamma}&=\gD\quad\text{on }\Gamma,
				   &
		u|_{\Gamma^c} &= 0\quad\text{on }\Gamma^c;\label{eqn:min_PDEb}
	\end{align}
	respectively.
	The dual norms~\eqref{eqn:H12_dual_norm_def} and~\eqref{eqn:H12_dual_alt_def}
	satisfy for any $\gN\in H^{-1/2}(\Gamma)$ that
	\begin{align}\label{eqn:Norm_HN_char}
		\NormHN[\Gamma]{\Omega}{\wn}{\gN}
		&=
		\max_{\substack{0\ne v\in H^1(\Omega)}}
		\frac{\left|\left\langle \gN,v\right\rangle_{\Gamma}\right|}	
			{\|v\|_{H^1(\Omega),\wn}}
		=\min_{\substack{\mathbf{m}\in H(\Omega,\Div) \\
		\mathbf{m}\cdot \nu=\gN\text{ on }\Gamma\\\mathbf{m}\cdot \nu=0\text{ on }\Gamma^c}}
		\|\mathbf{m}\|_{H(\Omega,\Div),\wn},\\
		\NormHNalt[\Gamma]{\Omega}{\wn}{\gN}
		&=
		\max_{\substack{0\ne v\in H^1(\Omega)\\v|_{\Gamma^c}=0}}
		\frac{\left|\left\langle \gN,v\right\rangle_{\Gamma}\right|}	
			{\|v\|_{H^1(\Omega),\wn}}
		=\min_{\substack{\mathbf{m}\in H(\Omega,\Div) \\
		\mathbf{m}\cdot \nu=\gN\text{ on }\Gamma}}
		\|\mathbf{m}\|_{H(\Omega,\Div),\wn}.\label{eqn:Norm_HN_alt_char}
	\end{align}
\end{theorem}
\begin{proof}
	\W{Given $\gD\in H^{1/2}(\Gamma)$, 
		the weak form of~\eqref{eqn:min_PDE} is the Euler-Lagrange equation for the 
		minimisation of the strictly convex and coercive norm $\|\bullet\|_{H^1(\Omega),\wn}$ 
		over the affine space $\{v\in H^1(\Omega)\ :\ v|_\Gamma=\gD\}$.
		Hence its unique solution realises the minimum in~\eqref{eqn:H12_norm_def}.}
		\W{On the other hand, for any $\gN\in H^{-1/2}(\Gamma)$,} 
		the unique solution $v\in H^1(\Omega)$ to
	\begin{align*}
		-\Delta v + \wn^2 v&=0\quad\text{in }\Omega,
						   &
		(\nabla v\cdot\nu)|_{\Gamma}&=\gN\quad\text{on }\Gamma,
				   &
		(\nabla v\cdot\nu)|_{\Gamma^c} &= 0\quad\text{on }\Gamma^c
	\end{align*}
	\W{realises (by~\eqref{eqn:Gauss_dualtiy}) the maximum in~\eqref{eqn:Norm_HN_char} with 
	$\NormHN[\Gamma]{\Omega}{\wn}\gN^2=\|v\|_{H^1(\Omega),\wn}^2
	=\left\langle \gN,\conj{v}\right\rangle_{\Gamma}$. 
	A similar argument identifies $\mathbf{m}\coloneqq\nabla v\in H(\Omega,\Div)$ 
	with $\|\mathbf{m}\|_{H(\Omega,\Div),\wn}=\|v\|_{H^1(\Omega),\wn}$ as the minimizer for the
	right-hand side of~\eqref{eqn:Norm_HN_char}. }
	The corresponding identities for the alternative norms~\eqref{eqn:H12_alt_def}--\eqref{eqn:H12_dual_alt_def}
	follow verbatim; further details are omitted.
\end{proof}
The point of the alternative norms~\eqref{eqn:H12_alt_def}--\eqref{eqn:H12_dual_alt_def}
is the duality of the boundary conditions on $\Gamma^c$ 
(cf.~\eqref{eqn:H12_norm_def}--\eqref{eqn:H12_dual_norm_def}
and~\eqref{eqn:Norm_HN_char}--\eqref{eqn:Norm_HN_alt_char}) 
as a result of~\eqref{eqn:Gauss_dualtiy}. 

The well-known Fourier transform characterisation of the trace space $H^{1/2}(\R^{n-1})$
provides an explicit characterisation of the 
minimal extension norm~\eqref{eqn:H12_norm_def} 
for the boundary $\partial \R^\dim_+=\R^{n-1}\times\{0\}$ (identified by $\R^{n-1}$)
of the half-space 
$\R^\dim_+\coloneqq\{(\overline x,x_n)\in\mathbb{R}^{n-1}\times\R\ :\ x_n>0\}$.
Let 
$\mathcal{F}\gD$ denote the $n-1$ dimensional Fourier transform of 
$\gD\in H^{1/2}(\R^{n-1})$. 

\begin{lemma}[Fourier representation]\label{lem:Fourier_halfspace}
	Any $\gD\in H^{1/2}(\R^{n-1})$ and $\sigma>0$ satisfy
	\begin{align*}
		\NormHD[\R^{n-1}]{\R^{\dim}_+}{\wn}\gD^2
		&=(2\pi)^{\dim-1}
		\int_{\R^{\dim-1}} (\wn^2 + |\xi|^2)^{1/2}|\mathcal{F}\gD|^2(\xi)\d \xi\\
		&\approx |\gD|_{H^{1/2}(\R^{n-1})}^2 + \wn\|\gD\|_{L^2(\R^{n-1})}^2.
	\end{align*}
	The constants hidden in the notation ``$\approx$'' exclusively depend on~$n$.
\end{lemma}
\begin{proof}
	Observe that \cref{lem:minimal extension norm} holds verbatim 
	for the half-space $\R^{n}_+$ with (non-compact) boundary $\R^{n-1}$.
	Lemma 2.7.1 in~\cite{Say:RetardedPotentialsTime2016}
	establishes the identity
	\begin{align*}
		\|u\|_{H^1(\R^{\dim}_+),\wn}^2 
		=(2\pi)^{\dim-1}\int_{\R^{\dim-1}} 
			(\wnAbs^2 + |\xi|^2)^{1/2}|\mathcal{F}\gD|^2(\xi)\d \xi
	\end{align*}
	for the solution $u\in H^1(\R^{\dim}_+)$ to~\eqref{eqn:min_PDE} with 
	$\gD\in H^{1/2}(\Gamma)$ for $\Gamma=\R^{n-1}$.
	Since $u$ minimises the norm 
	$\|u\|_{H^1(\R^\dim_+),\wn}=\NormHD[\R^{\dim-1}]{\R^\dim_+}{\wn}\gD$ 
	by \cref{lem:minimal extension norm},
	this proves the first equality.
	The well-known Plancherel identity~\cite[thm.~8.29]{Fol:RealAnalysisModern1999} and the
	characterisation ~\cite[prop.~3.4]{DPV:HitchhikersGuideFractional2012} of the Gagliardo seminorm for $\R^{n-1}$ reveal
	\begin{align*}
		\|\gD\|_{L^2(\R^{n-1})}^2 
		= \int_{\R^{\dim-1}} |\mathcal{F}\gD|^2(\xi)\d \xi
		\quad\text{and}\quad
		|\gD|_{H^{1/2}(\R^{n-1})}^2
		&= C(n) \int_{\R^{\dim-1}} |\xi||\mathcal{F}\gD|^2(\xi)\d \xi
	\end{align*}
	for an explicitly known constant $C(n)$ that exclusively depends on $n$.
	This and $\tfrac12(\wnAbs+|\xi|)^2\leq \wnAbs^2+|\xi|^2 \leq (\wnAbs+|\xi|)^2$
	(from a Cauchy inequality) imply the final equivalence.
\end{proof}

The following key result compares the minimal extension norm for connected extension sets
with the intrinsic weighted \W{norm $\NormGD[\Gamma]{\wn}\bullet$ from~\eqref{eqn:Sob_norm_def} that
generalises the classical Sobolev-Slobodeckij norm 
$\NormHclassic[\Gamma]\bullet=\NormGD[\Gamma]{1}\bullet$.}
The connectedness assumption poses no loss of generality as discussed in \cref{rem:multiply_connected} below: 
minimal extension norms naturally decompose as a sum over the connectivity components,
and the analysis for connected sets extends directly to the general case.%
\begin{theorem}[relation with intrinsic norm]\label{thm:characterisation_intrinsic}
	Let the Lipschitz set $\Omega\subset\R^n$ be connected and 
	$\Gamma\subset\partial\Omega$ an isolated boundary component.
	Any $\gD\in H^{1/2}(\Gamma)$ satisfies
	\begin{enumerate}[label=(\roman*)]
		\item $\displaystyle
			\Const{rel}^{-1}\NormGD[\Gamma]{\wn}\gD
		\leq \NormHD[\Gamma]{\Omega}{\wn}\gD
		\leq\Const{rel}\NormGD[\Gamma]{\max\{1,\wn\}}\gD
		$
	\end{enumerate}
	and, if additionally $\partial\Omega\setminus\Gamma\ne\emptyset$, then
	\begin{enumerate}[label=(\roman*)]\setcounter{enumi}{1}
		\item $\displaystyle
			\Const{rel}^{-1}\NormGD[\Gamma]{\max\{1,\wn\}}\gD
		\leq \NormHDalt[\Gamma]{\Omega}{\wn}\gD
		\leq\Const{rel}\NormGD[\Gamma]{\max\{1,\wn\}}\gD.
		$
	\end{enumerate}
	The constant $\Const{rel}>0$ is independent of $\wn>0$ and exclusively depends on $\Omega$ and $\Gamma$.
\end{theorem}
\begin{proof}%
	The proof considers any $\gD\in H^{1/2}(\Gamma)$ and splits into 3 steps.

	\noindent\emph{Step 1} provides 
	the lower bound in \emph{(i)}, namely
	\begin{align}\label{eqn:char_step_1}
		\NormGD[\Gamma]{\wn}\gD
		\leq \Const{rel}\NormHD[\Gamma]{\Omega}{\wn}\gD.
	\end{align}
	Observe that $\Omega$ may be replaced by $B_R\cap \Omega$
	for some ball $B_R=B_R^n$ with sufficiently large radius $R>0$ 
	such that $B_R\cap \Omega$ is connected and $\Gamma\subset \partial(B_R\cap\Omega)$
	is an isolated boundary component of $B_R\cap \Omega$.
	It then follows from~\eqref{eqn:H12_norm_def} that
	\begin{align}
		\NormHD[\Gamma]{B_R\cap \Omega}{\wn}\gD
		\leq 
		\NormHD[\Gamma]{\Omega}{\wn}\gD.
	\end{align}
	Hence the remaining parts of this step assume $\Omega$ to be bounded.
	Consider the minimal extension $u\in H^1(\Omega)$ 
	of $\gD=u|_{\Gamma}$ with $\|u\|_{H^1(\Omega),\wn}=\NormHD[\Gamma]{\Omega}{\wn}\gD$
	from \cref{lem:minimal extension norm}.

	\medskip
	\noindent\emph{Case 1a} ($\wn\leq 1$).
	Let $\langle u \rangle\coloneqq\frac{1}{|\Omega|}\int_{\Omega}u\d x$ denote
	the integral mean of $u$.
	The scaling of the weighted Sobolev-Slobodeckij norm for $\wn\leq 1$ and the 
	trace inequality~\eqref{eqn:C_tr_1} provide
	\begin{align*}
		\NormGD[\Gamma]{\wn}{\gD-\langle u \rangle}
		\leq\NormGD[\partial\Omega]{1}{\gD-\langle u \rangle}
		\leq \Const{tr,1}\|u-\langle u \rangle\|_{H^1(\Omega),1}
		\leq\Const{tr,1}\Const{P,\Omega}\|\nabla u\|_{L^2(\Omega)}
	\end{align*}
	with a Poincar\'e inequality~\eqref{eqn:C_P} in the last step.
	Since $[\langle u \rangle]_{H^{1/2}(\Gamma)}=0$ vanishes 
	for the constant $\langle u \rangle$, the
	definition~\eqref{eqn:Sob_norm_def} of $\NormGD[\Gamma]{\wn}\bullet$ for $\wn\leq 1$
	reveals
	\begin{align*}
		\NormGD[\Gamma]{\wn}{\langle u \rangle}
		=\wn\,|\Gamma|^{1/2}|\langle u \rangle|
		=\wn\,\frac{|\Gamma|^{1/2}}{|\Omega|^{1/2}}\|\langle u \rangle\|_{L^2(\Omega)}.
	\end{align*}
	A triangle inequality and the stability
	$\|\langle u \rangle\|_{L^2(\Omega)}\leq \|u\|_{L^2(\Omega)}$
	of the integral mean leads with the two previous
	displayed formulae to
	\begin{align*}
		\NormGD[\Gamma]{\wn}{\gD}
		\leq\Const{tr,1}\Const{P,\Omega}\|\nabla u\|_{L^2(\Omega)}
			+\frac{|\Gamma|^{1/2}}{|\Omega|^{1/2}}\, \wn\|u\|_{L^2(\Omega)}
		\leq\Const{rel}\|u\|_{H^1(\Omega),\wn}
	\end{align*}
	for $\Const{rel}^2=\Const{tr,1}^2\Const{P,\Omega}^2+ |\Gamma| / |\Omega|$ 
	from a Cauchy inequality in the last step.
	This and $\|u\|_{H^1(\Omega),\wn}=\NormHD[\Gamma]{\Omega}{\wn}\gD$
	prove~\eqref{eqn:char_step_1} for any $\wn\leq 1$.

	\medskip
	\noindent\emph{Case 1b} ($\wn\geq 1$).
	The multiplicative trace inequality for the weighted norms in%
	~\cite[lem.~3.1 and cor.~3.2]{MS:ConvergenceAnalysisFinite2010}
	provide (for $k_0=1$ therein) 
	a constant $\Const{tr,MS}$ that exclusively depends on $\Omega$ with
	\begin{align*}
		\NormGD[\partial\Omega]{ \wn}\gD
		\leq\Const{tr,MS}\|u\|_{H^1(\Omega),\wn}
		\qquad\text{for all }\wn\geq1.
	\end{align*}
	This and $\NormGD[\Gamma]{ \wn}\bullet\leq \NormGD[\partial\Omega]{ \wn}\bullet$
	verify~\eqref{eqn:char_step_1} with $\Const{rel}=\Const{tr,MS}$
	for all $\wn\geq1$ and conclude step 1.

	\medskip
	\noindent\emph{Step 2} establishes the first inequality in \emph{(ii)} 
	for $\Gamma^c=\partial\Omega\setminus\Gamma\ne\emptyset$, namely
	\begin{align}\label{eqn:char_step_2}
		\NormGD[\Gamma]{\max\{1,\wn\}}\gD
		\leq \Const{rel}\NormHDalt[\Gamma]{\Omega}{\wn}\gD.
	\end{align}
	The definitions~\eqref{eqn:H12_norm_def} and~\eqref{eqn:H12_alt_def} of 
	the minimal extension norms as infima reveal
	\begin{align}\label{eqn:Norm_vs_alt}
		\NormHD[\Gamma]{\Omega}{\wn}\gD\leq
		\NormHDalt[\Gamma]{\Omega}{\wn}\gD.
	\end{align}
	Hence~\eqref{eqn:char_step_2} follows for $\wn\geq1$ from step 1.
	To prove~\eqref{eqn:char_step_2} for $\wn\leq 1$,
	let $u\in H^1(\Omega)$ denote the minimal extension of $\gD=u|_{\Gamma}$
	with $u|_{\Gamma^c}=0$ on $\Gamma^c=\partial\Omega\setminus\Gamma$ and
	$\|u\|_{H^1(\Omega),\wn}=\NormHDalt[\Gamma]{\Omega}{\wn}\gD$
	from \cref{lem:minimal extension norm}.
	For a ball $B_R=B_R^n$ with sufficiently large radius $R>0$, 
	the bounded subset $U\coloneqq B_R\cap \Omega$ is connected and
	its boundary contains the isolated components 
	$\Gamma, \Gamma^c\subset\partial(B_R\cap \Omega)$.
	The unweighted trace inequality~\eqref{eqn:C_tr_1} reveals
	\begin{align*}
		\Const{tr,1}^{-2}\NormHclassic[\Gamma]\gD^2
		\leq \|u\|_{H^1(U),1}^2 = \|\nabla u\|_{L^2(U)}^2 + \|u\|_{L^2(U)}^2.
	\end{align*}
	Since $U$ is bounded and $\Gamma^c\ne\emptyset$, 
	the Friedrichs inequality~\eqref{eqn:C_F} provides
	\begin{align*}
		\|u\|_{L^2(U)}\leq \Const{F}\|\nabla u\|_{L^2(U)}.
	\end{align*}
	The combination of the two previous displayed inequalities 
	with $\Const{rel}^2 = \Const{tr,1}^{2}(1+\Const{F}^2)$~reads
	\begin{align*}
		\Const{rel}^{-1}\NormHclassic[\Gamma]\gD
		\leq\|\nabla u\|_{L^2(U)}\leq\|u\|_{H^1(U),\wn}.
	\end{align*}
	This and 
	$\|u\|_{H^1(U),\wn}\leq\NormHD[\Gamma]{\Omega}{\wn}\gD$ 
	reveal~\eqref{eqn:char_step_2} for $\wn\leq 1$ and
	conclude step~2.

	\medskip
	\noindent\emph{Step 3} considers the remaining estimate
	\begin{align}\label{eqn:char_step_3}
		\NormHDalt[\Gamma]{\Omega}{\wn}\gD
		\leq\Const{rel}\NormGD[\Gamma]{\max\{1,\wn\}}\gD,
	\end{align}
	which in combination with~\eqref{eqn:Norm_vs_alt} provides the upper bound
	in \emph{(i)} and \emph{(ii)}.
	It suffices to prove~\eqref{eqn:char_step_3} for all $\wn\geq1$,
	since 
	$\NormHDalt[\Gamma]{\Omega}{\wn}\bullet\leq\NormHDalt[\Gamma]{\Omega}{1}\bullet$ 
	holds for all $\wn\leq 1$.
	Hence the remaining parts of this proof assume $\wn\geq1$.

	The key idea is a local flattening of the boundary 
	(e.g., as in the proof of~\cite[lem.~2.7.3]{Say:RetardedPotentialsTime2016})
	and the norm representation \cref{lem:Fourier_halfspace} for the half-space.
	By the Lipschitz property of $\Omega$, there exists $J\in\mathbb{N}$ open sets
	$U_1,\dots,U_J\subset\R^n$ covering $\Gamma\subset\partial\Omega$ 
	and bi-Lipschitz homeomorphisms $\Psi_j:B_1^n\to U_j$
	(both $\Psi_j$ and $\Psi_j^{-1}$ are Lipschitz)
	from the $n$-dimensional unit ball $B_1^n$ (about the origin) onto $U_j$ with
	\begin{align*}
		\Psi_j(\R^n_+\cap B_1^n) = \Omega\cap U_j\eqqcolon \Omega_j
		\qquad\text{and}\qquad
		\Psi_j(B_1^{n-1}) = \Gamma\cap U_j\eqqcolon \Gamma_j
	\end{align*}
	for all $j=1,\dots,J$ with the identification $B_1^{n-1}\cong\partial\R^n_+\cap B_1^n$.
	The restrictions of $\Psi_j$ to bi-Lipschitz homeomorphisms from
	$\R^n_+\cap B_1^n$ and $B_1^{n-1}$ onto $\Omega_j$ and $\Gamma_j$ are denoted by
	\begin{align*}
		\Psi_j^+:\R^n_+\cap B_1^n \to \Omega_j
		\qquad\text{and}\qquad
		\Psi_j^{\partial}:B_1^{n-1}\to \Gamma_j,
		\qquad\text{respectivley}.
	\end{align*}
	Consider smooth functions $\varphi_j\in C^\infty_0(U_j)$ that induce a
	partition of unity for $\Gamma_j$ for $j=1,\dots,J$,
	i.e., $\sum^{J}_{j=1} \varphi_j\equiv 1$ on $\Gamma$.
	Observe that $\varphi_j\circ \Psi_j\in C^\infty_0(B_1)$ is compactly supported in some smaller $n$-dimensional ball
	$B_r\subset B_1$ with radius $r<1$ for all $j=1,\dots,J$.
	Hence there exists $\psi\in C^\infty_0(B_1)$ with $\psi\equiv 1$ on $B_r$
	and compact support in $B_1$.

	The function
	$\gD^{(j)}\coloneqq (\varphi_j\gD)\circ \Psi_j^{\partial}\in H^{1/2}(B_1^{n-1})$ 
	is compactly supported in $B_{r}^{n-1}$ and extends
	by zero to $H^{1/2}(\R^{n-1})$ for any $j=1,\dots,J$.
	Let $u_j\in H^1(\R^n_+)$ satisfy 
	\begin{align}\label{eqn:uj_def}
		u_j|_{\R^{n-1}}=\gD^{(j)}
		\qquad\text{and}\qquad
		\|u_j\|_{H^1(\R^{n}_+),\wn} 
		=\NormHD[\R^{n-1}]{\R^{n}_+}{\wn}{\gD^{(j)}}
	\end{align}
	(its existence follows as in \cref{lem:minimal extension norm}).
	Since $\psi$ has compact support in $B_1^n$, $(\psi u_j)\circ (\Psi_j^+)^{-1}\in H^{1}(\Omega_j)$
	extends by zero to $H^1(\Omega)$.
	Moreover, $\psi\equiv 1$ on the support of $\varphi_j\circ \Psi_j$ implies
	$(\psi u_j)|_{\Gamma}=\gD^{(j)}$
	such that $v\coloneqq \sum^{J}_{j=1} (\psi u_j)\circ (\Psi_j^+)^{-1}$ satisfies
	$v|_{\Gamma} = \gD$ by the partition of unity property of $(\varphi_j)_{j=1}^J$.
	Hence~\eqref{eqn:H12_norm_def}, a triangle inequality, and a change of variables 
	with the bi-Lipschitz $\Psi_j$ show
	\begin{align}\label{eqn:gD_ref_bound}
		\NormHD[\Gamma]{\Omega}{\wn}\gD
		\leq\|v\|_{H^1(\Omega),\wn}
		\leq \sum^{J}_{j=1} \|(\psi u_j)\circ (\Psi_j^+)^{-1}\|_{H^1(U_j),\wn}
		\lesssim\sum^{J}_{j=1} \|\psi u_j\|_{H^1(\R^n_+),\wn}.
	\end{align}
	For any $j=1,\dots,J$, the product rule provides
	\begin{align}\label{eqn:phi_uj_bound}
		\|\psi u_j\|_{H^1(\R^n_+),\wn}
		&\leq \|\psi\|_{L^\infty(\R^n_+)}\|u_j\|_{H^1(\R^n_+),\wn}
		+\|\nabla \psi\|_{L^\infty(\R^n_+)}\|u_j\|_{L^2(\R^n_+)}\\
		&\leq \left(\|\psi\|_{L^\infty(\R^n_+)} + \wn^{-1}\|\nabla \psi\|_{L^\infty(\R^n_+)}\right)\nonumber
		\|u_j\|_{H^1(\R^n_+),\wn}
		\lesssim \NormHD[\R^{n-1}]{\R^n_+}{\wn}{\gD^{(j)}}
	\end{align}
	with~\eqref{eqn:uj_def} and $\wn^{-1}\leq 1$ from $1\leq \wn$ 
	by assumption in the last step.
	\Cref{lem:Fourier_halfspace} controls
	\begin{align}
		\NormHD[\R^{n-1}]{\R^{n}_+}{\wn}{\gD^{(j)}}^2
		\lesssim
		\left|\gD^{(j)}\right|_{H^{1/2}(\R^{n-1})}^2
		+\wn\left\|\gD^{(j)}\right\|_{L^2(B_1^{n-1})}^2
		.\label{eqn:M_to_Sob}
	\end{align}
	Since $\gD^{(j)}$ is
	supported in $B_r^{n-1}$, the Gagliardo seminorm splits~as
	\begin{align}%
		\left|\gD^{(j)}\right|_{H^{1/2}(\R^{n-1})}^2
		&= \left|\gD^{(j)}\right|_{H^{1/2}(B_1^{n-1})}^2
		+2\int_{B_r^{n-1}} \int_{\R^{n-1}\setminus B_1^{n-1}}
		\frac{\left|\gD^{(j)}(x)\right|^2}{|x-y|^{n}}\d y\d x,
		\label{eqn:Sob_split}
	\end{align}
	where $\left|\bullet\right|_{H^{1/2}(B_1^{n-1})}$ is defined as 
	in~\eqref{eqn:Gagliardo_def} for the boundary part $B_1^{n-1}$.
	An integration in spherical coordinates and $0<1-r\leq|x-y|$ 
	for any $x\in B_r^{n-1}$ and $y\in \R^{n-1}\setminus B_1^{n-1}$ show
	\begin{align*}
		\int_{\R^{n-1}\setminus B_1^{n-1}}
		{|x-y|}^{-n}\d y\d x
		\leq\int_{1-r}^\infty\int_{\partial B_1^{n-1}}\rho^{-n}\d\Theta\d \rho
		\leq |\partial B_1^{n-1}|\int_{1-r}^\infty\rho^{-2}\d \rho
		<\infty.
	\end{align*}
	Hence~\eqref{eqn:Sob_split} and $1\leq\wn$ result in
	\begin{align*}
		\left|\gD^{(j)}\right|_{H^{1/2}(\R^{n-1})}^2
		\lesssim \left|\gD^{(j)}\right|_{H^{1/2}(B_1^{n-1})}^2
		+\wn\left\|\gD^{(j)}\right\|_{L^2(B_1^{n-1})}^2.
	\end{align*}
	This leads in combination with~\eqref{eqn:M_to_Sob} and a change of
	variables for $\gD^{(j)}=(\varphi_j\gD)\circ \Psi_j^{\partial}$
	with the bi-Lipschitz $\Psi_j^{\partial}$ to
	\begin{align}\label{eqn:gD_Sob_bound}
		\NormHD[\R^{n-1}]{\R^{n}_+}{\wn}{\gD^{(j)}}
		\lesssim
		\left|\gD^{(j)}\right|_{H^{1/2}(B_1^{n-1})}^2
		+\wn\left\|\gD^{(j)}\right\|_{L^2(B_1^{n-1})}^2
		\lesssim
		\NormGD[\Gamma]{\wn}{\varphi_j\gD}.
	\end{align}
	It is well-know 
	(cf.~\cite{Tri:MultiplicationPropertiesSpaces1977,Zol:MultiplicationDansEspaces1977} for references) that the multiplication by a 
	smooth function defines a bounded operator in 
	$H^{1/2}(\Gamma)$ and in $L^2(\Gamma)$ such that 
	\begin{align}\label{eqn:multiply_smooth}
		\NormGD[\Gamma]{\wn}{\varphi_j\gD}
		\lesssim
		\NormHclassic[\Gamma]{\gD} + \wn\|\gD\|_{L^2(\Gamma)}
		\lesssim
		\NormGD[\Gamma]{\wn}{\gD}
	\end{align}
	with $1\leq \wn$ in the last step.
	The combination of~\eqref{eqn:gD_ref_bound}--\eqref{eqn:phi_uj_bound} 
	and~\eqref{eqn:gD_Sob_bound}--\eqref{eqn:multiply_smooth}
	establish~\eqref{eqn:char_step_3}
	and conclude the proof.
\end{proof}
	\Cref{thm:characterisation_intrinsic} enables the characterisation of~\eqref{eqn:H12_norm_def}
	for domains postulated in part \emph{(a)} of the introduction.
	The result for exterior domains employs a novel Friedrichs-type inequality for 
	$\operatorname{DtN}$ boundary conditions for
	$n\geq3$ from~\cite{GS:DirichlettoNeumannOperatorHelmholtz2025}.
	Recall that a domain is in particular connected.
\begin{theorem}[characterisation]\label{thm:characterisation}
	Any $\gD\in H^{1/2}(\Gamma)$ satisfies
	\begin{enumerate}[label=(\roman*)]
		\item if $\Omega$ is a bounded domain:
		\begin{align*}%
			\Const{char}^{-1}\NormGD[\Gamma]{\wn}\gD
			\leq \NormHD[\Gamma]{\Omega}{\wn}\gD
			\leq\Const{char}
			\NormGD[\Gamma]{\wn}\gD,
		\end{align*}
	\item if $\Omega$ is an exterior domain and $n\geq3$:
		\begin{align*}%
			\Const{char}^{-1}\NormGD[\Gamma]{\max\{1,\wn\}}\gD
			\leq \NormHD[\Gamma]{\Omega}{\wn}\gD
			\leq\Const{char}
			\NormGD[\Gamma]{\max\{1,\wn\}}\gD.
		\end{align*}
	\end{enumerate}
	The constant $\Const{char}>0$ is independent of $\wn>0$ and exclusively depends on $\Omega$ and $\Gamma$.
\end{theorem}
\begin{proof}%
	Consider any $\gD\in H^{1/2}(\Gamma)$.

	\medskip
	\noindent\emph{Part (i)} for bounded $\Omega$: \Cref{thm:characterisation_intrinsic}.i
	provides the lower bound as well as, in the case $\wn\geq1$, the upper bound
	for $\Const{char}=\Const{rel}$.
	Hence it remains to prove
	\begin{align}\label{eqn:char_part_i}
		\NormHD[\Gamma]{\Omega}{\wn}\gD
		\leq\Const{char}\NormGD[\Gamma]{\wn}\gD
	\end{align}
	for any $\wn\leq 1$.
	Consider the integral mean
	$\langle \gD \rangle\coloneqq\frac{1}{|\Gamma|}\int_\Gamma \gD\d s$.
	On one hand, the definition of the minimal extension norm~\eqref{eqn:H12_norm_def}
	with $\wn\leq 1$
	and the continuous right-inverse $R_\Gamma$ of the trace map onto $\Gamma$ 
	from~\eqref{eqn:R_Gamma_def} provide
	\begin{align*}
		\NormHD[\Gamma]{\Omega}{\wn}{\gD-\langle \gD \rangle}
		\leq \|R_{\Gamma}(\gD-\langle \gD \rangle)\|_{H^1(\Omega),1}
		\leq \Const{R,1}
		\NormHclassic[\Gamma]{\gD-\langle \gD \rangle}.
	\end{align*}
	Since $|\langle \gD \rangle|_{H^{1/2}(\Gamma)}=0$, 
	the Poincar\'e inequality~\eqref{eqn:Poincare_Gamma_def}
	on $\Gamma$ controls
	\begin{align*}
		\NormHclassic[\Gamma]{\gD-\langle \gD \rangle}
		\leq (1+\Const{P,\Gamma}^2)^{1/2}|\gD|_{H^{1/2}(\Gamma)}.
	\end{align*}
	On the other hand, the constant extension of $\langle \gD \rangle$ 
	by itself and~\eqref{eqn:H12_norm_def} show
	\begin{align*}
		\NormHD[\Gamma]{\Omega}{\wn}{\langle \gD \rangle}
		\leq \|\langle \gD \rangle\|_{H^{1}(\Omega,\wn)}
		=|\Omega|^{1/2} \,\wn |\langle \gD \rangle|
		=\frac{|\Omega|^{1/2}}{|\Gamma|^{1/2}}\,\wn\|\langle \gD \rangle\|_{L^2(\Gamma)}.
	\end{align*}
	The combination of the three previous displayed estimates, a triangle inequality,
	and the stability $\|\langle \gD \rangle\|_{L^2(\Gamma)}\leq\|\gD\|_{L^2(\Gamma)}$
	of the integral mean on $\Gamma$ reveal
	\begin{align*}
		\NormHD[\Gamma]{\Omega}{\wn}\gD
		\leq \Const{R,1}(1+\Const{P,\Gamma}^2)^{1/2}|\gD|_{H^{1/2}(\Gamma)}
		+ \frac{|\Omega|^{1/2}}{|\Gamma|^{1/2}}\,\wn\|\gD\|_{L^2(\Gamma)}.
	\end{align*}
	This and a \W{Cauchy inequality} verify~\eqref{eqn:char_part_i} for 
	$\Const{char}^2
	=\Const{R,1}^2(1+\Const{P,\Gamma}^2) + |\Omega|/|\Gamma|$
	and conclude the proof of \emph{(i)}.

	\medskip
	\noindent\emph{Part (ii)} for an exterior domain $\Omega$ and $n\geq3$:
	Since \cref{thm:characterisation_intrinsic} provides the upper bound and, 
	in the case $\wn\geq1$,
	the lower bound for $\Const{char}=\Const{rel}$, it remains to show
	\begin{align}\label{eqn:char_part_ii}
		\NormHclassic[\Gamma]\gD\leq \Const{char}\NormHD[\Gamma]{\Omega}{\wn}\gD
	\end{align}
	for $\wn\leq 1$.
	Let $u\in H^1(\Omega)$ denote the minimal extension of $\gD=u|_\Gamma$
	with $\NormHD[\Gamma]{\Omega}{\wn}\gD = \|u\|_{H^1(\Omega),\wn}$ from
	\cref{lem:minimal extension norm}.
	Since the (closed) complement $\R^n\setminus\Omega$ of the exterior domain $\Omega$ is bounded, 
	$\R^n\setminus\Omega\subset B_R$ holds for some ball $B_R$ with sufficiently large radius $R>0$.
	Then $U\coloneqq B_R\cap \Omega$ satisfies $\Gamma\subset\partial U$ and the
	unweighted trace inequality~\eqref{eqn:C_tr_1} provides
	\begin{align}\label{eqn:char_1_trace_inequality}
		\NormHclassic[\Gamma]\gD\leq\Const{tr,1}\|u\|_{H^1(U),1}.
	\end{align}
	By~\eqref{eqn:min_PDE} in \cref{lem:minimal extension norm}, the restriction 
	$u|_{B_R^+}$ to the (unbounded) exterior 
	$B_R^+\coloneq\R^n\setminus\overline{B_R}\subset \Omega$ of $B_R$
	coincides with the unique solution in $u^+\in H^1(B_R^+)$ to
	\begin{align}\label{eqn:exterior_problem}
		-\Delta u^+ + \wn^2 u^+ =0
		\quad\text{in }B_R^+
		\qquad\text{and}\qquad
		u^+|_{S_R} = h
		\quad\text{on }S_R=\partial B_R
	\end{align}
	for $h=u|_{S_R}\in H^{1/2}(S_R)$.
	The DtN operator $\operatorname{DtN}:H^{1/2}(S_R)\to H^{-1/2}(S_R)$ 
	from%
	~\cite[subsec.~3.2]{GS:DirichlettoNeumannOperatorHelmholtz2025}
	associated to this exterior Helmholtz problem maps any $h\in H^{1/2}(S_R)$
	to the normal trace of the unique solution to~\eqref{eqn:exterior_problem},
	i.e., $\operatorname{DtN} (u|_{S_R}) = \partial_r u|_{S_R}$ 
	for the normal (radial) derivative $\partial_r$ on $S_R$.
	Since $-\Delta u + \wn^2 u = 0$ on $U\subset\Omega$ by~\eqref{eqn:min_PDE} in 
	\cref{lem:minimal extension norm}, the Friedrichs-type inequality from%
	~\cite[eqn.~(3.7) in cor.~3.6]{GS:DirichlettoNeumannOperatorHelmholtz2025}
	provides a constant $\Const{F,DtN}>0$ that exclusively depends on $U$ and $R$ with
	\begin{align}\label{eqn:C_F_DtN}
		\|u\|_{L^2(U)}\leq \Const{F,DtN}\|u\|_{H^1(U),\wn}.
	\end{align}
	The combination of~\eqref{eqn:char_1_trace_inequality} and~\eqref{eqn:C_F_DtN}
	reveal for $\Const{char}\coloneqq \Const{tr,1}(1+\Const{F,DtN}^2)^{1/2}$
	that
	\begin{align*}
		\NormHclassic[\Gamma]\gD^2
		\leq \Const{tr,1}^2\|\nabla u\|_{L^2(U)}^2
			+\Const{tr,1}^2\Const{F,DtN}^2\|u\|_{H^1(U),\wn}^2
		\leq \Const{char}^2\|u\|_{H^1(U),\wn}^2.
	\end{align*}
	This and \W{$\|u\|_{H^1(U),\wn}\leq\|u\|_{H^1(\Omega),\wn}=\NormHD[\Gamma]{\Omega}{\wn}\gD$}
	verify~\eqref{eqn:char_part_ii} and conclude the proof.
\end{proof}

\begin{remark}[multiply connected Lipschitz sets]\label{rem:multiply_connected}
	\Cref{thm:characterisation_intrinsic,thm:characterisation}
	also characterise the minimal extension norms 
	for isolated boundary components $\Gamma\subset\partial\Omega$ of
	general (possibly multiply connected) Lipschitz sets $\Omega\subset\R^n$.
	Indeed, the minimal extension norms decompose as a sum
	over the $J\in \mathbb{N}$ connectivity components $\Omega_j$ 
	of $\Omega$ with $\Gamma\cap \partial\Omega_j\ne\emptyset$ (for $j=1,\dots,J$) in the sense that
	any $\gD\in H^{1/2}(\Gamma)$ satisfies%
	\begin{align*}
		\NormHD[\Gamma]{\Omega}{\wn}\gD^2
		=\sum^{J}_{j=1} \NormHD[\Gamma\cap \partial\Omega_j]{\Omega_j}{\wn}\gD^2
		\quad\text{and}\quad
		\NormHDalt[\Gamma]{\Omega}{\wn}\gD^2
		=\sum^{J}_{j=1} \NormHDalt[\Gamma\cap \partial\Omega_j]{\Omega_j}{\wn}\gD^2.
	\end{align*}
	For example, \cref{thm:characterisation}.i and~\eqref{eqn:Gagliardo_split} provide for bounded $\Omega$ and
	any $\gD\in H^{1/2}(\Gamma)$ that%
	\begin{align}\label{eqn:characterisation_bounded_multiply_connected}
		\min\{1,\wn\}^2\NormGD[\Gamma]{\wn}\gD^2
		\lesssim
		\NormHD[\Gamma]{\Omega}{\wn}\gD^2
		\approx \sum^{J}_{j=1} \NormGD[\Gamma\cap \partial\Omega_j]{\wn}\gD^2
		\leq \NormGD[\Gamma]{\wn}\gD^2.
	\end{align}
\end{remark}
\section[Wavenumber-explicit scaling and trace inequality]{$\wn$-explicit scaling and trace inequality}%
\label{sec:Properties of minimal extension norms}
This section discusses the scaling of the weighted minimal extension norms 
\W{$\NormH{\bullet}_{H^{\pm1/2}(\Gamma),\wn,\Omega}$ 
from~\eqref{eqn:H12_norm_def}--\eqref{eqn:H12_dual_norm_def}}
in the weight $\wn>0$
and establishes the $\wn$-explicit trace estimate.
Recall the lower bound $\wnlow=\wnLow^{-1}$ from~\eqref{eqn:wnlow_def} and
that $\Omega\subset\R^n$ denotes a Lipschitz set with 
isolated boundary component $\Gamma\subset\partial\Omega$ and its complement $\Gamma^c=\partial\Omega\setminus\Gamma$. 
\begin{lemma}[scaling in the weight]\label{lem:scaling in the weight}
	Any $\gD\in H^{1/2}(\Gamma)$ and $\gN\in H^{-1/2}(\Gamma)$ satisfy
	\begin{enumerate}[label=(\roman*)]
		\item $\displaystyle\begin{aligned}[t]
			\wnlow\NormHD[\Gamma]{\Omega}{1}\gD
			&\leq
			\NormHD[\Gamma]{\Omega}{\wn}\gD
			\leq \max\{1,\min\{\wn,\Const{sc}\wnAbs^{1/2}\}\}\NormHD[\Gamma]{\Omega}{1}\gD,\\
			\wnlow\NormHN[\Gamma]{\Omega}{\wn}\gN
			&\leq
			\NormHN[\Gamma]{\Omega}{1}\gN
			\leq\max\{1,\min\{\wn,\Const{sc}\wnAbs^{1/2}\}\}\NormHN[\Gamma]{\Omega}{\wn}\gN.
		\end{aligned}$\\[0.7em]
	and the analogical scaling holds for 
	$\NormHalt{\bullet}_{H^{\pm1/2}(\Gamma),\wn,\Omega}$
	from~\eqref{eqn:H12_alt_def}--\eqref{eqn:H12_dual_alt_def}.
	\item If $\Omega$ is an exterior domain and $n\geq3$, then
		\begin{align*}
			\NormHD[\Gamma]{\Omega}{1}\gD
			&\leq
			\Const{sc}\NormHD[\Gamma]{\Omega}{\wn}\gD
			\quad\text{and}\quad
			\NormHN[\Gamma]{\Omega}{\wn}\gN
			\leq
			\Const{sc}\NormHN[\Gamma]{\Omega}{1}\gN.
		\end{align*}
		\item If $\partial\omega\setminus\Gamma\ne\emptyset$ holds
			for every connectivity component $\omega$ of $\Omega$, then
		\begin{align*}
			\NormHDalt[\Gamma]{\Omega}{1}\gD
			&\leq
			\Const{sc}\NormHDalt[\Gamma]{\Omega}{\wn}\gD
			\quad\text{and}\quad
			\NormHNalt[\Gamma]{\Omega}{\wn}\gN
			\leq
			\Const{sc}\NormHNalt[\Gamma]{\Omega}{1}\gN.
		\end{align*}
	\end{enumerate}
	The constant $\Const{sc}>0$ is independent of $\wn>0$ and exclusively depends on $\Gamma$ and $\Omega$.
\end{lemma}
\begin{proof}
	The weighted $H^1$ norm~\eqref{eqn:H1_norm} scales in the weight $\wn>0$ as
	\begin{align*}
		\min\{1,\wnAbs\}\|v\|_{H^1(\Omega),1}
		\leq \|v\|_{H^1(\Omega),\wn}
		\leq\max\{1, \wnAbs\}\|v\|_{H^1(\Omega),1}
		\qquad\text{for all }v\in H^1(\Omega)
	\end{align*}
	and implies the analogical scaling for the minimal extension norm~\eqref{eqn:H12_norm_def}, namely
	\begin{align}\label{eqn:trivial_scaling_a}
		\min\{1,\wnAbs\}\NormHD[\Gamma]{\Omega}{1}\bullet
		&\leq \NormHD[\Gamma]{\Omega}{\wn}\bullet
		\leq\max\{1, \wnAbs\}\NormHD[\Gamma]{\Omega}{1}\bullet.
	\end{align}
	For large weights $1\leq\wn$, the second estimate in~\eqref{eqn:trivial_scaling_a} improves
	with \cref{thm:characterisation_intrinsic}.i as
	\begin{align*}
		\Const{rel}^{-1}\NormHD[\Gamma]{\Omega}{\wn}\bullet
		\leq \NormGD[\Gamma]{\wn}\bullet
		\leq \wn^{1/2}\NormHclassic[\Gamma]\bullet
		\leq \Const{rel}\wn^{1/2}\NormHD[\Gamma]{\Omega}{1}\bullet.
	\end{align*}
	By \cref{rem:operator_norm_est}, 
	this and~\eqref{eqn:trivial_scaling_a} verify \W{\emph{(i)}}
	for $\Const{sc}=\Const{rel}^2$. 

	The characterisations from \cref{thm:characterisation}.ii 
	and \cref{thm:characterisation_intrinsic}.ii 
	(with \cref{rem:operator_norm_est}) lead to \W{\emph{(ii)}--\emph{(iii)}}; 
	further details are omitted.
\end{proof}%
\W{
\begin{cor}[scaling in the weight]\label{cor:scaling in the weight}
	The trace norms
	satisfy~\eqref{eqn:scaling_HN_bounded}--\eqref{eqn:scaling_HN_bounded_improved}.
\end{cor}
\begin{proof}
	This is an immediate consequence of \cref{lem:scaling in the weight}
	and the well-known equivalence 
	$\NormHclassic[\Gamma]{\bullet}\approx \NormHD[\Gamma]{\Omega}{1}{\bullet}$ (cf.%
	~\cref{thm:characterisation_intrinsic}.i) and similarly 
	for the operator norms.
\end{proof}}

The remaining parts of this section are devoted to the $\wn$-explicit trace inequality in \cref{lem:trace_inequality}
below and depart with two preliminary results.
\begin{lemma}[comparison]\label{lem:norm characterisation}
	For any $\gD\in H^{1/2}(\Gamma)$, 
	the trace norms are equivalent with
	\begin{align}
		\NormHD[\Gamma]{\Omega}{\wn}{\gD}
		&\leq
		\NormHDalt[\Gamma]{\Omega}{\wn}{\gD}
		\leq
		\Const{cmp}
		\NormHD[\Gamma]{\Omega}{\max\{1,\wn\}}{\gD}
		\leq
		\Const{cmp}\wnLow
		\NormHD[\Gamma]{\Omega}{\wn}{\gD}\label{eqn:Norm_HD_alt_equiv}
	\end{align}
	If $\Omega$ is a bounded Lipschitz domain and $\Gamma^c\ne\emptyset$, 
	the constant function $\gD\equiv1$ %
	satisfies
	\begin{align}\label{eqn:comparison_sharp_a}
		\NormHD[\Gamma]{\Omega}{\wn}{1}
		&\leq \Const{rev}\wnlow\,
		\NormHDalt[\Gamma]{\Omega}{\wn}{1}.%
	\end{align}
	The constants $\Const{cmp},\Const{rev}>0$ are independent of $\wn>0$ and
	exclusively depend on $\Gamma$ and~$\Omega$.
\end{lemma}
\begin{proof}

	The first inequality in~\eqref{eqn:Norm_HD_alt_equiv} is clear from the definition of the norms as infima.
	Consider any $\gD\in H^{1/2}(\Gamma)$ and any connectivity component $\omega$ of $\Omega$.
	If $\partial\omega\setminus\Gamma=\emptyset$, 
	the definition in~\eqref{eqn:H12_norm_def} and~\eqref{eqn:H12_alt_def} shows
	\begin{align}\label{eqn:H_empty_equal}
		\NormHDalt[\Gamma\cap \partial\omega]{\omega}{\wn}\gD
		=
		\NormHD[\Gamma\cap \partial\omega]{\omega}{\wn}\gD.
	\end{align}
	In the case $\partial\omega\setminus\Gamma\ne\emptyset$, the characterisation 
	\cref{thm:characterisation_intrinsic} 
	provides
	\begin{align*}
		\Const{rel}^{-1}\NormHDalt[\Gamma\cap \partial\omega]{\omega}{\wn}\gD
		&\leq
		\NormGD[\Gamma\cap \partial\omega]{\max\{1,\wn\}}\gD
		\leq
		\Const{rel}\NormHD[\Gamma\cap \partial\omega]{\omega}{\max\{1,\wn\}}\gD.
	\end{align*}
	Since the minimal extension norms decompose as a sum over the 
	connectivity components of the extension set
	(cf.%
	~\cref{rem:multiply_connected} for details), 
	this %
	establishes the second inequality
	in~\eqref{eqn:Norm_HD_alt_equiv} for 
	$\Const{cmp}\coloneqq\max\{1,\Const{rel}^2\}$.
	The final inequality is follows from~\eqref{eqn:scaling_HN_bounded} 
	\W{(cf.~\cref{cor:scaling in the weight})}.
	
	Assume that $\Omega\subset\R^n$ is a bounded domain with $\Gamma^c\ne\emptyset$.
	On one hand, the constant function is an admissible extension of its trace in $H^1(\Omega)$ 
	such that~\eqref{eqn:H12_norm_def} reveals
	\begin{align}\label{eqn:const_function}
		\NormHD[\Gamma]{\Omega}{\wn}{1}\leq \|1\|_{H^1(\Omega),\wn} = \wnAbs\,\|1\|_{L^2(\Omega)}
		=\wnAbs\,|\Omega|^{1/2}.
	\end{align}
	On the other hand, \cref{thm:characterisation_intrinsic}.ii, $|1|_{H^{1/2}(\Gamma)}=0$,
	and \cref{rem:multiply_connected} provide
	\begin{align*}
		\max\{1,\wn\}|\Gamma|^{1/2}
		= \NormGD[\Gamma]{\max\{1,\wn\}}1
		\leq \Const{rel}\NormHD[\Gamma]{\Omega}{\max\{1,\wn\}}1.
	\end{align*}
	This and~\eqref{eqn:const_function} verify~\eqref{eqn:comparison_sharp_a}
	for $\Const{rev}^2=\Const{rel}^2|\Omega|/|\Gamma|$ and conclude the proof.
\end{proof}

\begin{theorem}[different extension sets]\label{thm:extension_sets}
	Let $\Omega,\widetilde{\Omega}\subset\R^n$ be Lipschitz sets and $\Gamma\subset\partial\Omega\cap \partial\widetilde{\Omega}$
	an isolated boundary component of $\partial\Omega$ and of $\partial\widetilde{\Omega}$.
	Any $\gD\in H^{1/2}(\Gamma)$ satisfies
	\begin{align}\label{eqn:norm_comparison}
		\NormHDalt[\Gamma]{\widetilde{\Omega}}{\wn}\gD
		&\leq \Const{eq}\,\NormHD[\Gamma]{\Omega}{\max\{1,\wn\}}\gD
		\leq \Const{eq}\wnLow\,\NormHD[\Gamma]{\Omega}{\wn}\gD.
	\end{align}
	If additionally
	\begin{enumerate}[label=(\roman*)]
		\item $\Omega$ is connected and $\widetilde{\Omega}$ is bounded, then
		\begin{align*}%
			\NormHD[\Gamma]{\widetilde{\Omega}}{\wn}\gD
			&\leq \Const{eq}\NormHD[\Gamma]{\Omega}{\wn}\gD,
		\end{align*}
		\item $\Omega$ is an exterior domain and $n\geq3$, then
		\begin{align*}%
			\NormHDalt[\Gamma]{\widetilde{\Omega}}{\wn}\gD
			&\leq \Const{eq}\NormHD[\Gamma]{\Omega}{\wn}\gD,
		\end{align*}
		\item $\Omega$ is connected and $\partial\Omega\setminus\Gamma\ne\emptyset$, then
		\begin{align*}%
			\NormHDalt[\Gamma]{\widetilde{\Omega}}{\wn}\gD
			&\leq \Const{eq}\NormHDalt[\Gamma]{\Omega}{\wn}\gD.
		\end{align*}
	\end{enumerate}
	The constant $\Const{eq}>0$ is independent of $\wn>0$ and exclusively depends on 
	$\Gamma$, $\Omega$, and $\widetilde{\Omega}$.
\end{theorem}
\begin{proof}
	Recall that the norms $\NormHD[\Gamma]{\widetilde{\Omega}}{\wn}\bullet\leq\NormHDalt[\Gamma]{\widetilde{\Omega}}{\wn}\bullet$ 
	(by~\eqref{eqn:Norm_HD_alt_equiv}) decompose as a sum over the connectivity components of $\widetilde{\Omega}$ by
	\cref{rem:multiply_connected} and coincide 
	on the connectivity components $\omega$ with $\partial\omega\setminus\Gamma=\emptyset$ 
	by~\eqref{eqn:H_empty_equal}.
	Hence \cref{thm:characterisation_intrinsic} and~\eqref{eqn:Gagliardo_split} lead to
	\begin{align}\label{eqn:extension_sets_A}
		\NormHD[\Gamma]{\widetilde{\Omega}}{\wn}\gD
		\leq
		\NormHDalt[\Gamma]{\widetilde{\Omega}}{\wn}\gD
		\leq \Const{rel}\NormGD[\Gamma]{\max\{1,\wn\}}\gD.
	\end{align}
	Similarly, \cref{thm:characterisation_intrinsic}.i and~\eqref{eqn:Gagliardo_split} reveal for the $J\in\mathbb{N}$
	connectivity components $\Omega_j$ of $\Omega$ with $\Gamma\cap\Omega_j\ne\emptyset$ for $j=1,\dots,J$ that
	\begin{align}\nonumber
		\max\{1,\Const{spl}\}^{-1}&\NormGD[\Gamma]{\max\{1,\wn\}}\gD
		\leq
		\sum^{J}_{j=1} \NormGD[\Gamma\cap \partial\Omega_j]{\max\{1,\wn\}}\gD\\
		&\leq 
		\Const{rel} \NormHD[\Gamma]{\Omega}{\wn}\gD
		\leq \wnLow\Const{rel}
		\NormHD[\Gamma]{\Omega}{\wn}\gD\label{eqn:extension_sets_B}
	\end{align}
	with~\eqref{eqn:scaling_HN_bounded} \W{from \cref{cor:scaling in the weight}}
	in the last step.
	The combination of~\eqref{eqn:extension_sets_A}--\eqref{eqn:extension_sets_B} verify~\eqref{eqn:norm_comparison}
	for $\Const{eq}=\max\{1,\Const{spl}\}\Const{rel}^2$.
	Assume for the remaining parts of the proof that $\Omega$ is connected.
	For bounded $\widetilde{\Omega}$,~\emph{(i)} follows with~\eqref{eqn:characterisation_bounded_multiply_connected}
	from \cref{thm:characterisation}.i and \cref{thm:characterisation_intrinsic}.i for
	$\Const{eq}=\Const{char}\Const{rel}$.
	If $\widetilde{\Omega}$ is an exterior domain, \eqref{eqn:extension_sets_A} and \cref{thm:characterisation}.ii 
	reveal~\emph{(ii)}.
	Finally, if \W{$\Gamma^c\ne\emptyset$}, 
	\eqref{eqn:extension_sets_A} 
	and \cref{thm:characterisation_intrinsic}.ii establish~\emph{(iii)}
	and conclude the proof.
\end{proof}

	Recall the Sobolev spaces $H^1_{0,\gamma}(\Omega)$ and $H_{0,\gamma}(\Omega,\Div)$
	with homogeneous boundary conditions on $\gamma\subset\partial\Omega$ from
	\cref{sub:General notation}.
	\begin{theorem}[$\wn$-explicit trace inequality]\label{lem:trace_inequality}
	Let $\Omega,\widetilde{\Omega}\subset\R^n$ be Lipschitz sets and $\Gamma\subset\partial\Omega\cap \partial\widetilde{\Omega}$
	with oriented unit normal $\nu$
	be an isolated boundary component of $\partial\Omega$ and of $\partial\widetilde{\Omega}$.
	Any $v\in H^1(\Omega\cup\widetilde{\Omega})$ and
	$\mathbf{m}\in H(\Omega\cup\widetilde{\Omega},\Div)$ satisfy
		\begin{enumerate}[label=(\roman*)]
			\item 
			$\displaystyle 
			\NormHD[\Gamma]{\Omega}{\wn}{v|_{\Gamma}}\leq \|v\|_{H^1(\Omega),\wn}$,
			\item[(ii.a)] 
			$\displaystyle 
			\NormHD[\Gamma]{\Omega}{\wn}{{ v}|_{\Gamma}}
			\leq \Const{tr}\,\|{ v}\|_{H^1(\widetilde{\Omega}),\max\{1,\wn\}}
			\leq \Const{tr}\wnLow\,\|{ v}\|_{H^1(\widetilde{\Omega}),\wn},$
			\item[(ii.b)] if $\widetilde{\Omega}$ is connected and $\Omega$ is bounded or \\
					if $\widetilde{\Omega}$ is an exterior domain and $n\geq3$ or\\
					if $\widetilde{\Omega}$ is connected, $\partial\widetilde{\Omega}\setminus\Gamma\ne\emptyset$, and 
					$ v\in H_{0,\partial\widetilde{\Omega}\setminus\Gamma}(\widetilde{\Omega})$, then\\
			$\displaystyle 
			\NormHD[\Gamma]{\Omega}{\wn}{{ v}|_{\Gamma}}
			\leq \Const{tr}\,\|{ v}\|_{H^1(\widetilde{\Omega}),\wn},$
		\item[(iii.a)]
			$\displaystyle 
			\NormHN[\Gamma]{\Omega}{\wn}{(\mathbf{m}\cdot \nu)|_{\Gamma}}
			\leq \Const{tr}\wnLow\|\mathbf{m}\|_{H(\Omega,\Div),\wn}
			$,
		\item[(iii.b)]
				if $\mathbf{m}\in H_{0,\partial\Omega\setminus\Gamma}(\Omega,\Div)$, then\\
			$\displaystyle 
			\NormHN[\Gamma]{\Omega}{\wn}{(\mathbf{m}\cdot \nu)|_{\Gamma}}
			\leq 
				\|\mathbf{m}\|_{H(\Omega,\Div),\wn}
				$,
		\item[(iv.a)]
			$\displaystyle 
			\NormHN[\Gamma]{\Omega}{\wn}{({{\mathbf{m}}}\cdot \nu)|_{\Gamma}}
			\leq \Const{tr}\wnLow\,\|{{\mathbf{m}}}\|_{H(\widetilde{\Omega},\Div),\wn}$,
		\item[(iv.b)] if $\widetilde{\Omega}$ is bounded, $\Omega$ is connected, and 
			$\mathbf{m}\in H_{0,\partial\widetilde{\Omega}\setminus\Gamma}(\widetilde{\Omega},\Div)$ or\\
					if $\Omega$ is an exterior domain and $n\geq3$, then\\
			$\displaystyle 
			\NormHN[\Gamma]{\Omega}{\wn}{({{\mathbf{m}}}\cdot \nu)|_{\Gamma}}
			\leq \Const{tr}\,\|{{\mathbf{m}}}\|_{H(\widetilde{\Omega},\Div),\wn}$.
		\end{enumerate}
		The constant $\Const{tr}>0$ is independent of $\wn>0$ and 
		exclusively depends on $\Omega$ and $\widetilde{\Omega}$.
	\end{theorem}
	\begin{proof}
		A straightforward combination of \cref{thm:extension_sets} 
		with the definitions~\eqref{eqn:H12_norm_def}
		and~\eqref{eqn:H12_alt_def} of $\NormHD[\Gamma]{\omega}{\wn}\bullet\leq\NormHDalt[\Gamma]{\omega}{\wn}\bullet$
		as infima over (a subset of) $H^1(\omega)$ for $\omega\in\{\Omega,\widetilde{\Omega}\}$
		lead~to~\emph{(i)}--\emph{(ii)} for $\Const{tr}=\Const{eq}$.

		To prove~\emph{(iii)}--\emph{(iv)}, consider any $\mathbf{m}\in H(\Omega\cup\widetilde{\Omega},\Div)$ and
		set $\gN\coloneqq (\mathbf{m}\cdot\nu)|_{\Gamma}\in H^{-1/2}(\Gamma)$.
		By \cref{rem:operator_norm_est}, \cref{lem:norm characterisation,thm:extension_sets} reveal
		\begin{enumerate}[align=left,leftmargin=*]
			\item[ad \emph{(iii.a)}:] $\displaystyle
				\NormHN[\Gamma]{\Omega}{\wn}\gN\leq\Const{cmp}\wnLow\NormHNalt[\Gamma]{\Omega}{\wn}\gN$,
			\item[ad \emph{(iv.a)}:] $\displaystyle
				\NormHN[\Gamma]{\Omega}{\wn}\gN\leq\Const{eq}\wnLow\NormHNalt[\Gamma]{\widetilde{\Omega}}{\wn}\gN$
		\end{enumerate}
		and, under the additional assumptions on $\Omega$ and $\widetilde{\Omega}$ in \emph{(iv.b)},
		\begin{enumerate}[align=left,leftmargin=*]
			\item[ad \emph{(iv.b)}] if $\widetilde{\Omega}$ is bounded and $\Omega$ is connected: $\displaystyle
				\NormHN[\Gamma]{\Omega}{\wn}\gN\leq\Const{eq}\NormHN[\Gamma]{\widetilde{\Omega}}{\wn}\gN$,
			\item[ad \emph{(iv.b)}] if $\Omega$ is an exterior domain and $n\geq3$: $\displaystyle
				\NormHN[\Gamma]{\Omega}{\wn}\gN\leq\Const{eq}\NormHNalt[\Gamma]{\widetilde{\Omega}}{\wn}\gN$.
		\end{enumerate}
		This and the identification of the operator norms as minima over 
		(a subset of) $H(\omega,\Div)$ for
		$\omega\in\{\Omega,\widetilde{\Omega}\}$ in \cref{lem:minimal extension norm} result in \emph{(iii)}--\emph{(iv)}
		and conclude the proof for $\Const{tr}\coloneqq\max\{\Const{eq},\Const{cmp}\}$.
	\end{proof}
	Transmission or interface problems are often posed in $\R^n\setminus\Gamma$
	for an interface $\Gamma=\partial\Omega_{}$ that is boundary of some 
	bounded Lipschitz domain $\Omega\subset\R^n$ with connected exterior.
	In this case, the trace inequality in \cref{lem:trace_inequality} simplifies
	for the extension set $\Omega$.
	\begin{cor}[$\wn$-explicit trace inequality for interface problems]\label{cor:trace inequality}
		Let $\Omega\subset\R^n$ be a bounded Lipschitz set with unit normal $\nu$ on $\Gamma\coloneqq\partial\Omega$
		and connected exterior $\Omega_{\mathrm{ext}}\coloneqq\R^n\setminus\overline{\Omega}$.
		Any $v\in H^1(\R^n\setminus\Gamma)$ and $\mathbf{m}\in H(\R^n\setminus\Gamma,\Div)$ satisfy
		with $\Const{tr}>0$ from \cref{lem:trace_inequality} that
		\begin{enumerate}[label=(\roman*)]
			\item 
				$\displaystyle 
				\NormHD[\Gamma]{\Omega}{\wn}{v|_{\Gamma}}\leq \min\{\|v\|_{H^1(\Omega),\wn},
				\Const{tr}\|v\|_{H^1(\Omega_{\mathrm{ext}}),\wn}\}
				$,
			\W{\item 
				$\displaystyle 
				\NormHD[\Gamma]{\Omega_{\mathrm{ext}}}{\wn}{v|_{\Gamma}}
				\leq \min\{\Const{tr}\wnLow\|v\|_{H^1(\Omega),\wn},\|v\|_{H^1(\Omega_{\mathrm{ext}}),\wn} \}
		$,}
			\item $\displaystyle
				\NormHN[\Gamma]{\Omega}{\wn}{(\mathbf{m}\cdot \nu)|_{\Gamma}}
				\leq
				\min\{\|\mathbf{m}\|_{H(\Omega,\Div),\wn},\Const{tr}\wnLow\|\mathbf{m}\|_{H(\Omega_{\mathrm{ext}},\Div),\wn}\}$\W{,
			\item $\displaystyle
				\NormHN[\Gamma]{\Omega_{\mathrm{ext}}}{\wn}{(\mathbf{m}\cdot \nu)|_{\Gamma}}
				\leq
			\min\{\Const{tr}\|\mathbf{m}\|_{H(\Omega,\Div),\wn},\|\mathbf{m}\|_{H(\Omega_{\mathrm{ext}},\Div),\wn}\}$.}
		\end{enumerate}
	\end{cor}
\begin{proof}
	This is a direct consequence of \cref{lem:trace_inequality}; further details are omitted.
\end{proof}

\section{Application to boundary integral operators}%
\label{sub:Operator bounds for boundary integral operators}
This section considers some %
Lipschitz set $\Omega\subset\R^n$ with oriented unit normal $\nu$ on the boundary $\Gamma\coloneqq\partial\Omega$
and a wavenumber $s\in\mathbb{C}$ with positive real part $\Re s>0$.
Recall the Dirichlet and Neumann jumps and averages~\eqref{eqn:jump_mean_def}--\eqref{eqn:N_jump_average}.
The single and double layer potentials~\eqref{eqn:opSD} map $H^{\pm1/2}(\Gamma)$ into $H^1(\R^n\setminus\Gamma)$ and
have the integral representations%
\begin{align}
	(\opS g)(x)
	&\coloneqq\int_{\R^n\setminus\Gamma} g(y)\, G(x,y)\d y
	&&\text{for all }x\in\R^n\setminus\Gamma,\label{eqn:opS_def}\\
	(\opD g)(x)
	&\coloneqq\int_{\R^n\setminus\Gamma} g(y)\, \nabla_{y}G(x,y)\cdot \nu(y)\d y
	&&\text{for all }x\in\R^n\setminus\Gamma,\label{eqn:opD_def}
\end{align}
for sufficiently smooth $g:\Gamma\to\mathbb{C}$ with
the fundamental solution $G$ 
~\cite[chap.~9]{McL:StronglyEllipticSystems2000}
of the Helmholtz operator~\eqref{eqn:HH_op}, i.e., $\opL G(x,\bullet)=\delta_{x}$ for
the Dirac delta function $\delta_x$.
Their (averaged) traces define
the \emph{single layer}, \emph{double
layer}, \emph{dual double layer}, and \emph{hypersingular boundary integral operators}%
~\cite[eqn.~(3.6)]{SS:BoundaryElementMethods2011}
for any $\gD\in H^{1/2}(\Interface)$ and $\gN\in H^{-1/2}(\Interface)$~by
\begin{align}\label{eqn:Bop_V}
	\opV\phantom{'} &: H^{-1/2}(\Interface)\to H^{1/2}(\Interface)
		 &&\text{with}\qquad
		 \opV \gN\coloneqq\meanD[\Gamma]{\opS \gN},\\\label{eqn:Bop_K}
	\opK\phantom{'} &: H^{1/2}(\Interface)\to H^{1/2}(\Interface)
		 &&\text{with}\qquad
		 \opK \gD\coloneqq \meanD[\Gamma]{\opD \gD},\\\label{eqn:Bop_Kdual}
	\opKdual &: H^{-1/2}(\Interface)\to H^{-1/2}(\Interface)
			 &&\text{with}\qquad
			 \opKdual \gN\coloneqq \meanN[\Gamma]{\opS \gN},\\\label{eqn:Bop_W}
	\opW &: H^{1/2}(\Interface)\to H^{-1/2}(\Interface)
		 &&\text{with}\qquad
		 \opW \gD\coloneqq \meanN[\Gamma]{\opD \gD}.
\end{align}
It is well-known (for $\Re s>0$) that the sesquilinear form associated to $\opL$ from~\eqref{eqn:HH_op}
\begin{align*}
	\opl(v,w)\coloneqq\int_{\R^n\setminus\Gamma}\nabla v\cdot \nabla\conj w + s^2 v\, \conj w\d x
	\qquad\text{for all }v,w\in H^1(\R^n\setminus\Gamma)
\end{align*}
is coercive~\cite{BH:FormulationVariationnelleEspacetemps1986} (see also%
~\cite[lem.~3.2]{FHS:SkeletonIntegralEquations2024}) with 
\begin{align}\label{eqn:opl_coercive}
	|\opl(v,v)|\geq \Re\left(\frac{\conj s}{\W{\wn}}\opl(v,v)\right)\geq \CLinv \|v\|_{H^1(\R^n\setminus\Gamma),\W{\W{\wn}}}^2
	\qquad\text{for all }v\in H^1(\R^n\setminus\Gamma).
\end{align}
The single and double layer potentials~\eqref{eqn:opSD} with~\eqref{eqn:opS_def}--\eqref{eqn:opD_def}
as well as the boundary layer operators~\eqref{eqn:Bop_V}--\eqref{eqn:Bop_W} 
\W{admit explicit and uniform bounds in terms of the inverse coercivity constant $\W{\wn}/\Re s\geq1$}
with respect to the corresponding weighted norms.%
\begin{theorem}[$\wn$-explicit continuity estimates]\label{thm:opSD_bounds}
	Let $s\in\mathbb{C}$ with $\Re s>0$ and set $\wn=|s|$. 
	\W{$\opV:H^{-1/2}(\Omega)\to H^{1/2}(\Gamma)$
		and $\opW:H^{1/2}(\Omega)\to H^{-1/2}(\Gamma)$
	are isomorphisms with
	\begin{align}\label{eqn:V_coercive}
		\Re\left(\frac{\conj s}{\W{\wn}}\left\langle\opV\gN,\conj \gN \right\rangle_{\Gamma} \right)
		&\geq\frac{\Re s}{\W{\wn}}
		\NormHN[\Gamma]{\Omega}{\wn}{\gN}^2
		&&\text{for all }\gN\in H^{-1/2}(\Gamma),\\
		\Re\left(\frac{\conj s}{\W{\wn}}\left\langle\opW\gD,\conj \gD \right\rangle_{\Gamma} \right)
		&\geq\frac{\Re s}{\W{\wn}}
		\NormHD[\Gamma]{\Omega}{\wn}{\gD}^2
		&&\text{for all }\gD\in H^{1/2}(\Gamma).
	\end{align}
	Moreover, any $\gD\in H^{1/2}(\Gamma)$ and any $\gN\in H^{-1/2}(\Gamma)$ satisfy}
	\begin{align}\label{eqn:opS_bound}
		\|\opS \gN\|_{H^1(\R^n\setminus\Gamma),\wn}
		&\leq
		\CLInv\NormHN[\Gamma]{\Omega}{\wn}{\gN},\\
		\NormHD[\Gamma]{\Omega}{\wn}{\opV \gN}
		&\leq
		\CLInv\NormHN[\Gamma]{\Omega}{\wn}{\gN},\label{eqn:opV_bound}\\
		\NormHN[\Gamma]{\Omega}{\wn}{\opKdual \gN}
		&\leq
		\left(\frac{1}{2}+\CLInv\right)\NormHN[\Gamma]{\Omega}{\wn}{\gN}
		\label{eqn:opKdual_bound},\\
		\NormHD[\Gamma]{\Omega}{\wn}{\opW^{-1} \gN}
		&\leq
		\CLInv\NormHN[\Gamma]{\Omega}{\wn}{\gN},\label{eqn:opW_inv_bound}\\
		\|\opD \gD\|_{H^1(\R^n\setminus\Gamma),\wn}
		&\leq
		\CLInv\NormHD[\Gamma]{\Omega}{\wn}{\gD},\label{eqn:opD_bound}\\
		\NormHD[\Gamma]{\Omega}{\wn}{\opK \gD}
		&\leq
		\left(\frac{1}{2}+\CLInv\right)\NormHD[\Gamma]{\Omega}{\wn}{\gD},\label{eqn:opK_bound}\\
		\NormHN[\Gamma]{\Omega}{\wn}{\opW \gD}
		&\leq
		\CLInv\NormHD[\Gamma]{\Omega}{\wn}{\gD}\label{eqn:opW_bound},\\
		\NormHN[\Gamma]{\Omega}{\wn}{\opV^{-1} \gD}
		&\leq
		\CLInv\NormHD[\Gamma]{\Omega}{\wn}{\gD}.\label{eqn:opV_inv_bound}
	\end{align}
\end{theorem}
The proof of \cref{thm:opSD_bounds} below requires a well-known~\cite{BH:FormulationVariationnelleEspacetemps1986}
bound for the normal derivative of $\opL$-harmonic functions in the current setting.
\begin{lemma}[bound for the normal derivative]\label{lem:normal_bound}
	Let $s\in\mathbb{C}\setminus\{0\}$ and set $\wn=|s|$. 
	Any $v\in H^1(\Omega)$ with $-\Delta v+s^2 v=0$ in $\Omega$ satisfies
	\begin{align*}
		\NormHN[\Gamma]{\Omega}{\wn}{(\partial_{\nu}v)|_{\Gamma}}
		\leq \|v\|_{H^1(\Omega),\wn}.
	\end{align*}
\end{lemma}
\begin{proof}
	This is \cref{lem:trace_inequality}.iii.b applied to $\mathbf{m}\coloneqq\nabla v$
	with $\|\mathbf{m}\|_{H(\Div,\Omega),\wn}=\|v\|_{H^1(\Omega),\wn}$.
\end{proof}
\begin{proof}[Proof of \cref{thm:opSD_bounds}]
	This proof follows standard
	arguments~\cite{BH:FormulationVariationnelleEspacetemps1986,McL:StronglyEllipticSystems2000,Say:RetardedPotentialsTime2016} 
	(see also~\cite[appendix A]{FHS:SkeletonIntegralEquations2024}).
	The single and double layer potentials satisfy
	\begin{align*}
		\jumpD[\Gamma]{\opS\gN} 
		&= 0,
		&
		\jumpN[\Gamma]{\opS\gN} 
		&=\pm\gN&&\text{for all }\gN\in H^{-1/2}(\Gamma),\\
		\jumpD[\Gamma]{\opD\gD} 
		&= \pm\gD,
		&
		\jumpN[\Gamma]{\opD\gD} 
		&=0&&\text{for all }\gD\in H^{1/2}(\Gamma)
	\end{align*}
	by~\cite[thm.~6.11]{McL:StronglyEllipticSystems2000},
	where the sign $\pm$ depends on the orientation of the unit normal~$\nu$.
	A standard (piecewise) integration by parts for $v\coloneqq\opS\gN$ 
	in~\eqref{eqn:opl_coercive} leads 
	with $ v|_{\Omega}=v|_{\R^n\setminus\overline{\Omega}}$ by the jump conditions to
	\begin{align}\label{eqn:single_layer_proof}
		\CLinv \|v\|_{H^1(\R^n\setminus\Gamma),\wn}^2
		\leq \Re\left(\frac{\conj s}{\W{\wn}}\left\langle \jumpN[\Gamma] v,\conj v|_{\Gamma}
			\right\rangle_{\Gamma} \right)
			\leq \NormHN[\Gamma]{\Omega}{\wn}{\gN}\NormHD[\Gamma]{\Omega}{\wn}{v|_{\Gamma}}.
	\end{align}
	This and 
	$\NormHD[\Gamma]{\Omega}{\wn}{v|_{\Gamma}}\leq \|v\|_{H^1(\R^n\setminus\Gamma),\wn}$ from
	\cref{lem:trace_inequality}.i establish~\eqref{eqn:opS_bound}.
	\W{A trace inequality with \cref{lem:trace_inequality}.i 
	for $v|_{\Gamma}=\opV\gN$ and~\eqref{eqn:opS_bound} show~\eqref{eqn:opV_bound}.
	Since $\jumpN[\Gamma]v=\gN$ by the jump conditions,~\eqref{eqn:single_layer_proof}
	reveals the coercivity~\eqref{eqn:V_coercive} of the single layer operator
	and thus the bound~\eqref{eqn:opV_inv_bound} for its inverse $\opV^{-1}$.}
	Moreover,
	\begin{align*}
		\W{\opKdual\gN = (\partial_\nu\opS\gN|_{\Omega})|_{\Gamma} - 
			\frac{1}{2}\jumpN[\Gamma]{\opS\gN}
			= (\partial_\nu v|_{\Omega})|_{\Gamma} - \frac{1}{2}\gN}
	\end{align*}
	follows from the jump conditions and the properties of the jumps and averages.
	It is well-known~\cite{Say:RetardedPotentialsTime2016}
	that the single and double layer potentials 
	map into the $\opL$-harmonic functions, 
	\W{i.e.,} $-\Delta v + s^2 v = 0$ in $\R^n\setminus\Gamma$.
	Hence a triangle inequality, and \cref{lem:normal_bound} show
	\begin{align*}
		\NormHN[\Gamma]{\Omega}{\wn}{\opKdual\gN}
		\leq\frac{1}{2} \NormHN[\Gamma]{\Omega}{\wn}{\gN} + \|v\|_{H^1(\Omega),\wn}.
	\end{align*}
	This and~\eqref{eqn:opS_bound} verify~\eqref{eqn:opKdual_bound}.
	The proof of the remaining bounds~\eqref{eqn:opD_bound}--\eqref{eqn:opW_bound}
	related to the double layer potential \W{and its traces} follow verbatim; 
	further details are omitted.
\end{proof}

\W{
	\Cref{thm:characterisation_intrinsic,thm:characterisation}
	(and \cref{rem:multiply_connected})
	characterise the weighted trace norms through the intrinsic norm
	\eqref{eqn:Sob_norm_def}. Combined with \cref{thm:opSD_bounds}, this gives
	the $\sigma$-dependence of the potential and boundary integral operators, up to
	the coercivity factor $\Re s/\wn$ from \eqref{eqn:opl_coercive}.
	The corresponding bounds in the classical trace norms are known
	from~\cite[thms.~4.6--4.7]{BS:IntegralEquationMethods2022} and%
	~\cite[lem.~5.2]{FHS:SkeletonIntegralEquations2024}, 
	and are recovered from \cref{thm:opSD_bounds} 
	with the scaling relation~\eqref{eqn:scaling_HN_bounded}.
	For instance, the potential bounds from%
	~\cite{FHS:SkeletonIntegralEquations2024}
	read (in this notation) as
	\begin{align}\label{eqn:S_bound_classic}
		\|\opS \gN\|_{H^1(\R^n\setminus\Gamma),\wn}
		&\leq
		C\max\{1,\W{\wn}^{-1}\}\CLInv\NormHclassicDual[\Gamma]{\gN},\\
		\|\opD \gD\|_{H^1(\R^n\setminus\Gamma),\wn}
		&\leq
		C\max\{1,\W{\wn}^{1/2}\}\CLInv\NormHclassic[\Gamma]{\gD}.\label{eqn:D_bound_classic}
	\end{align}
	For $n\ge3$, the following corollary improves these bounds 
	in the low-frequency regime.
	Recall the (unweighted) classical Sobolev-Slobodeckij norms 
	$\NormH{\bullet}_{H^{\pm1/2}(\Gamma)}$ 
	for $H^{\pm1/2}(\Gamma)$.
\begin{cor}[unweighted trace norms]\label{thm:opSD_bounds_old}
	Assume $\Gamma=\partial\Omega$ is the boundary of an exterior Lipschitz domain 
	$\Omega\subset\R^n$ for $n\geq 3$, and
	let $s\in\mathbb{C}$ with $\Re s>0$ and $\wn=|s|$. 
	There exists $C$ exclusively depending on $\Gamma$ such that any $\gD\in H^{1/2}(\Gamma)$
	and $\gN\in H^{-1/2}(\Gamma)$ satisfy
\begingroup
\allowdisplaybreaks
	\begin{align*}
		\Re\left(\frac{\conj s}{\W{\wn}}\left\langle\opV\gN,\conj \gN \right\rangle_{\Gamma} \right)
		&\geq C\min\{1,\W{\wn}^{-1}\}\frac{\Re s}{\W{\wn}}
		\NormHclassicDual[\Gamma]{\gN}^2,\\
		\Re\left(\frac{\conj s}{\W{\wn}}\left\langle\opW\gD,\conj \gD \right\rangle_{\Gamma} \right)
		&\geq C\frac{\Re s}{\W{\wn}}
		\NormHclassic[\Gamma]{\gD}^2,\\
		\|\opS \gN\|_{H^1(\R^n\setminus\Gamma),\wn}
		&\leq
		C\CLInv\NormHclassicDual[\Gamma]{\gN},\\
		\NormHclassic[\Gamma]{\opV \gN}
		&\leq
		C\CLInv\NormHclassicDual[\Gamma]{\gN},\\
		\NormGN[\Gamma]{1}{\opKdual \gN}
		&\leq
		C\max\{1,\W{\wn}^{1/2}\}\frac{\W{\wn}}{\Re s}\NormHclassicDual[\Gamma]{\gN},\\
		\NormHclassic[\Gamma]{\opW^{-1} \gN}
		&\leq
		C\CLInv\NormHclassicDual[\Gamma]{\gN},\\
		\|\opD \gD\|_{H^1(\R^n\setminus\Gamma),\wn}
		&\leq
		C\max\{1,\W{\wn}^{1/2}\}\CLInv\NormHclassic[\Gamma]{\gD},\\
		\NormHclassic[\Gamma]{\opK \gD}
		&\leq
		C\max\{1,\W{\wn}^{1/2}\}\frac{\W{\wn}}{\Re s}\NormHclassic[\Gamma]{\gD},\\
		\NormGN[\Gamma]{1}{\opW \gD}
		&\leq
		C\max\{1,\W{\wn}\}\CLInv\NormHclassic[\Gamma]{\gD},\\
		\NormGN[\Gamma]{1}{\opV^{-1} \gD}
		&\leq
		C\max\{1,\W{\wn}\}\CLInv\NormHclassic[\Gamma]{\gD}.%
	\end{align*}
\endgroup
\end{cor}
\begin{proof}
	Observe that one has the freedom to choose the extension set $\Omega$ on 
	\enquote{either side} of the interface $\Gamma$ in \cref{thm:opSD_bounds}, 
	and here we consider the exterior domain $\Omega$.
	Since $n\geq3$, the combination of
	\cref{thm:opSD_bounds} and~\eqref{eqn:scaling_HN_bounded_improved}
	from \cref{cor:scaling in the weight} concludes the proof.
\end{proof}
}%
\begin{remark}[comparison]
	\W{The weighted trace norms are induced by the 
	Helmholtz energies~\eqref{eqn:H1_norm}--\eqref{eqn:Hdiv_norm}
	and measure boundary data at the natural PDE scale.
	This enables coercivity and continuity estimates 
	for the potential and boundary integral operators
	in \cref{thm:opSD_bounds} without the additional frequency factors present with
	the classical trace norms.
	They may therefore be viewed as energy norms for Dirichlet and
	Neumann data.

	In contrast, the classical 
	bounds~\cite{BH:ExistenceHmatrixApproximants2003,BS:IntegralEquationMethods2022}
	recovered with~\eqref{eqn:scaling_HN_bounded} may
	deteriorate as $\W{\wn}\to 0$ independently of the coercivity factor $\Re s/\W{\wn}$
	see, e.g.,~\eqref{eqn:S_bound_classic}.
	At least for $n\geq3$, \cref{thm:opSD_bounds_old} removes this critical low-frequency
	deterioration. The remaining factors of $\max\{1,\W{\wn}^{1/2}\}$ are exactly those dictated 
	by the scaling in \cref{lem:scaling in the weight}.

	These comparison factors may \emph{not} be sharp for slowly oscillatory boundary data. 
	Indeed, if $|\gD|_{H^{1/2}(\Gamma)}\lesssim \W{\wn}^\alpha\|\gD\|_{L^2(\Gamma)}$ 
	with $\alpha\in(0,1]$ and $\W{\wn}<1$,
	\cref{thm:opSD_bounds} for bounded $\Omega$ and the characterisation%
	~\eqref{eqn:Sob_norm_def} provide the sharper energy bound
	\begin{align*}
		\|\opD \gD\|_{H^1(\R^n\setminus\Gamma),\wn}
		&\lesssim \frac{\W{\wn}^{1+\alpha}}{\Re s}\|\gD\|_{L^2(\Gamma)}
	\end{align*}
	as opposed to the classical bound of the form
	$\lesssim\frac{\W{\wn}}{\Re s}\|\gD\|_{L^2(\Gamma)}$ from
	\cref{thm:opSD_bounds_old}.

	Extensions to heterogeneous, uniformly positive $L^\infty$
	coefficients} in the Helmholtz operator as in%
	~\cite{FHS:SkeletonIntegralEquations2024} are straightforward (with constants
	depending on the coefficients).
	\W{Analogous results} for the indefinite case $\Re s=0$ are stated in%
	~\cite[sec.~4]{GHS:StableSkeletonIntegral2025}.
\end{remark}

\section*{Acknowledgement}%
The author thanks Prof.~Stefan Sauter (Universität Zürich) for
insightful discussions and valuable suggestions that contributed to the development of this work.

	\bibliographystyle{alphaabbr}
	\bibliography{Bibliography}

\end{document}